\documentclass[a4paper]{article}
\usepackage{graphicx}
\usepackage{hyperref}
\usepackage{latexsym}
\newtheorem{theorem}{Theorem}
\newtheorem{corollary}{Corollary}
\begin{document}
\title{\LARGE{The significance of Nathanson's \textit{boss} factor in legitimising Aristotle's particularisation} \\ \vspace{+1ex} \large{Why we need to revise current interpretations of Cantor's, G\"{o}del's, Turing's and Tarski's formal reasoning}}
\author{\normalsize{Bhupinder Singh Anand}}
\date{\normalsize \textit {Draft of \today}}
\maketitle

\begin{abstract}
\footnotetext{Alix Comsi Internet Services Pvt.\ Ltd.\ Postal address: 32 Agarwal House, D Road, Churchgate, Mumbai - 400 020, Maharashtra, India.\ Email: re@alixcomsi.com, anandb@vsnl.com.}
\noindent I show---contrary to common beliefs tolerated by the ``bosses"---that any interpretation of ZF that admits Aristotle's particularisation is not sound; that the standard interpretation of PA is not sound; that PA is consistent but \(\omega\)-\textit{in}consistent; that a sound finitary interpretation of PA is definable in terms of Turing-computability; and that PA cannot be consistently extended to ZF.
\end{abstract}

\section{Preamble: The semantic and logical paradoxes}

We are all familiar with the semantic and logical paradoxes\footnote{Commonly referred to as the paradoxes of `self-reference', even though not all of them involve self-reference, e.g., the paradox constructed by Yablo \cite{Ya93}.} which involve---either implicitly or explicitly---quantification over an infinitude.

Where such quantification is not, or cannot be, explicitly defined in formal logical terms---eg. the classical expression of the Liar paradox as `This sentence is a lie'---the paradoxes per se cannot be considered as posing serious linguistic or philosophical concerns\footnote{It \textit{would} be a matter of serious concern if the word `This' in the English language sentence `This sentence is a lie' could be validly viewed as implicitly implying that: (i) there is a constructive infinite enumeration of English language sentences; (ii) to each of which a truth-value can be constructively assigned by the rules of a two-valued logic; and, (iii) in which `This' refers uniquely to a particular sentence in the enumeration. In \cite{Go31} G\"{o}del used the above perspective: (a) to show how the infinitude of formulas in a formally defined Peano Arithmetic P (\cite{Go31}, pp.9-13) could be constructively enumerated and referenced uniquely by natural numbers (\cite{Go31}, p.13-14); (b) to show how PA-provability values could be constructively assigned to PA-formulas by the rules of a two-valued logic (\cite{Go31}, p.13); and, (c) to construct a PA-formula that interprets as an arithmetical proposition that could be viewed as expressing the sentence `This P-sentence is P-unprovable' (cf.\ \cite{Go31}, p.37, footnote 67) without inviting a `Liar' type of contradiction.}, except to illustrate the absurd extent to which languages of common discourse need to tolerate ambiguity\footnote{In a lighter vein: Such absurdity is also highlighted by the universal appreciation of Charles Dickens' Mr. Bumble's retort that ``The law is an ass" - a quote oft used to refer to the absurdities that \href{http://www.shazbot.com/lawass/}{sometimes surface} in cases when judicial pronouncements attempt to resolve an ambiguity by subjective fiat that appeals to the powers - and duties - bestowed upon them for the practical resolution of precisely such an ambiguity, even when it may be theoretically irresolvable!}, both, for ease of expression, and for practical - even if not theoretically unambiguous and effective - communication in non-critical cases amongst intelligences capable of a lingua franca.

Now, addressing such ambiguity in critical cases\footnote{Such as communication between mechanical artefacts; or a putative communication between terrestrial and extra-terrestrial intelligences.} is the very raison d'etre of mathematical activity which is, first, the construction of richer and richer mathematical languages that can express those of our abstract concepts that can be subjectively\footnote{In what follows, I use the words `subjective' and `non-constructive' inter-changeably.} addressed unambiguously; and, thereafter, the study of the ability of the mathematical languages to precisely express and objectively\footnote{In what follows, I use the words `objective' and `constructive' inter-changeably.} communicate these concepts effectively.

However, even where the quantification can be made explicit---eg.
Russell's paradox\footnote{Define \(S\) by \{\(All \hspace{+.5ex} x: x \in S\) iff \(x \notin x\)\}; then \(S \in S\) iff \(S \notin S\).} or Yablo's paradox\footnote{Defining \(S_{i}\) for all \(i \geq 0\) as `For all \(j>i\), \(S_{j}\) is not true' seems to lead to a contradiction (cf.\ \cite{Ya93}).}---the question arises whether such quantification is constructive or not\footnote{For instance, in Russell's case it could be argued that the contradiction itself establishes that \(S\) cannot be constructively defined over the range of the quantifier. In Yablo's case it could be argued that truth values cannot be constructively assigned to any sentence covered by the quantification since, in order to decide whether \(S_{i}\) is true or not for any given \(i \geq 0\), we first need to decide whether \(S_{i+1}\) is true or not.}.

There are two issues involved here---not necessarily independent.

\subsection{Is quantification currently interpreted constructively over an infinite domain?}

The first---and more significant---is whether the currently accepted interpretations of formal quantification\footnote{Essentially as defined by Hilbert in terms of his \(\epsilon\)-function \cite{Hi27}} over an infinite domain can be treated as constructive or not.

Now, Brouwer\footnote{\cite{Br08}.} emphatically---and justifiably so far as number theory was concerned (see Appendix A)---asserted that Hilbert's interpretations of formal quantification were non-constructive. Although Hilbert's formalisation of the quantifiers (an integral part of his formalisation of Aristotle's logic of predicates) appeared adequate, Brouwer rejected Hilbert's interpretations of them on the grounds that the interpretations were open to ambiguity\footnote{They could not, therefore, be accepted as admitting effective communication.}.

However, Brouwer's objection was seen---also justifiably, as I show in a subsequent section---as unconvincingly rejecting a comfortable interpretation that---despite its Platonic overtones---appeared intuitively plausible to the larger body of academics which was increasingly attracted to, and influenced by, the remarkably expressive powers provided by Cantor-inspired set theories such as ZF.

Since Brouwer far preceded Turing, he was unable to offer his critics an alternative---and intuitively convincing---constructive definition of quantification.

Moreover, since Brouwer's objections did not gain much currency amongst mainstream logicians, they were unable to influence Turing, who---as I show in this paper---could easily have provided the necessary constructive interpretations sought by Hilbert for number theory, had he not himself been influenced by G\"{o}del's powerful presentation---and G\"{o}del's Platonic interpretation---of his own formal reasoning in G\"{o}del's seminal 1931 paper on formally undecidable arithmetical propositions\footnote{\cite{Go31}.}.

Thus, in his 1939 paper\footnote{\cite{Tu39}.} on ordinal-based logics, Turing applied the computational method that he described in his 1936 paper\footnote{\cite{Tu36}.} in seeking a categorical interpretation of Cantor's ordinal arithmetic (as defined in a set theory such as ZF)---rather than in seeking a categorical interpretation of PA---and apparently viewed his 1936 paper\footnote{\cite{Tu36}.} as complementing and extending G\"{o}del's and Cantor's reasoning.

Turing thus overlooked the fact that---as I show in a subsequent section---his 1936 paper\footnote{\cite{Tu36}.} actually conflicts with G\"{o}del's and Cantor's interpretations of their own, formal, reasoning by admitting an objective definition of satisfaction that yields a sound, finitary, interpretation \(\mathcal{I}_{PA(\beta rouwer/Turing)}\) of PA.

(In other words, whereas G\"{o}del's and Cantor's reasoning implicitly presumes that satisfaction under the standard interpretation \(\mathcal{I}_{PA(Standard/Tarski)}\) of PA can only be defined non-constructively in terms of subjectively verifiable truth, satisfaction under \(\mathcal{I}_{PA(\beta rouwer/Turing)}\) is defined constructively in terms of objectively verifiable Turing-computability.)

As a result, current theory continued---and continues to this day---to essentially follow Hilbert's Platonically-influenced definitions and interpretations of the quantifiers (based on Aristotle's logic of predicates) when defining them under the standard interpretation\footnote{cf.\ \cite{Me64}, p.107; \cite{Sh67}, p.23, p.209; \cite{BBJ03}, p.104. For purposes of this paper, I treat \cite{Me64} as a reliable and representative exposition---where cited---of the standard logical concepts and theory that are addressed at that point.} \(\mathcal{I}_{PA(Standard/Tarski)}\) of formal number theory.

The latter definitions and interpretations\footnote{eg.\ \cite{Me64}, pp.49-53.} are, in turn, founded upon Tarski's analysis of the inductive definability of the truth of compound expressions of a symbolic language under an interpretation in terms of the satisfaction of the atomic expressions of the language under the interpretation\footnote{\cite{Ta33}. Note that Tarski defines the formal sentence \(P\) as True if and only if \(p\)---where \(p\) is the proposition expressed by \(P\). In other words, the sentence ``Snow is white" is True if, and only if, it is \textit{subjectively} true in all cases; and it is \textit{subjectively} true in a particular case if, and only if, it expresses the \textit{subjectively verifiable} fact that snow \textit{is} white in that particular case. Thus, for Tarski the commonality of the satisfaction of the atomic formulas of a language under an interpretation (cf.\ \cite{Me64}, p.51(i)) is axiomatic.}.

\subsection{When is the concept of a completed infinity consistent with a formal language?}

The second issue is when, and whether, the concept of a completed infinity is consistent with the interpretation of a formal language.

Clearly, the consistency of the concept would follow immediately in any sound interpretation\footnote{We define an interpretation \(I_{S}\) of a formal system \(S\) as \textit{sound} if, and only if, every provable formula \([F]\) of \(S\) translates as a true proposition \(F\) under the interpretation (cf.\ \cite{BBJ03}, p.174).} of the axioms (and rules of inference) of a set theory such as ZF (whether such an interpretation exists at all is, of course, another question; one, moreover, that I address in a subsequent section).

In view of the perceived power of ZF\footnote{More accurately, ZFC.} as an unsurpassed language of rich and adequate expression of mathematically expressible abstract concepts precisely, it is not surprising that many of the semantic and logical paradoxes depend on the implicit assumption that the domain over which the paradox quantifies can always be treated as a well-defined mathematical object that can be formalised in ZF, even if this domain is not explicitly defined set-theoretically.

This assumption is rooted in the questionable belief that ZF can express all mathematical `truths'\footnote{In a subsequent section, I show how---particularly in the case of Goodstein's Theorem---such a belief leads to conclusions that fall short of accepted standards of mathematical rigour.}.

From this it is but a short step to the non-constructive argument---rooted in G\"{o}del's Platonic interpretation of his own formal reasoning in his 1931 paper\footnote{\cite{Go31}.}---that PA must have non-standard models.

In a series of recent papers\footnote{In particular \cite{An08f} and \cite{An08h}.} I have argued that both of the above foundational issues need to be reviewed carefully, and that we need to recognize explicitly the limitations on the ability of highly expressive mathematical languages such as ZF to communicate effectively, and the limitations on the ability of effectively communicating mathematical languages such as PA to adequately express abstract concepts---such as those involving Cantor's first limit ordinal \(\omega\).

Prima facie, the semantic and logical paradoxes seem to arise out of a blurring of this distinction, and an attempt to ask of a language more than it is designed to deliver.

\section{Introduction}

I now highlight five `beliefs' in the common foundations of philosophy, logic, mathematics and computability theory---resting upon Aristotle's logic of predicates and tolerated (even when not subscribed to) by the ``bosses"\footnote{The word ``bosses" is not intended to be pejorative, but to refer to the collective of reputed---and respected---experts in any field of human endeavour in the sense of Nathanson's comments in \cite{Na08}.}---which illustrate the point sought to be made by Melvyn B. Nathanson in his Opinion piece, ``Desperately Seeking Mathematical Truth", in the August 2008 Notices of the American Mathematical Society:

\begin{quote}
`` ... many great and important theorems don't actually have proofs. They have sketches of proofs, outlines of arguments, hints and intuitions that were obvious to the author (at least, at the time of writing) and that, hopefully, are understood and believed by some part of the mathematical community.

But the community itself is tiny. In most fields of mathematics there are few experts. Indeed, there are very few active research mathematicians in the world, and many important problems, so the ratio of the number of mathematicians to the number of problems is small. In every field, there are ``bosses" who proclaim the correctness or incorrectness of a new result, and its importance or unimportance.

Sometimes they disagree, like gang leaders fighting over turf. In any case, there is a web of semi-proved theorems throughout mathematics. Our knowledge of the truth of a theorem depends on the correctness of its proof and on the correctness of all of the theorems used in its proof. It is a shaky foundation."
\end{quote}

I shall show that standard arguments for the following are founded on plausible, but logically insecure, reasoning (and shall consider some curious consequences):

\begin{tabbing}
\noindent (a) \= A sound interpretation of ZF can appeal to Aristotle's particularisation; \\

\noindent (b) PA is \(\omega\)-consistent; \\

\noindent (c) The standard interpretation of first-order Peano Arithmetic (PA) is sound; \\

\noindent (d) PA is not categorical; and \\

\noindent (e) PA can be consistently extended to ZF.
\end{tabbing}

I shall show how, and why, such beliefs---instances of the influence of Nathan- son's ``boss" factor---are suspect for drawing sound conclusions in areas that cannot be expected to critically scrutinise their foundational framework.

For  instance, it is reasonable to assume that these foundational beliefs---which together imply that effective and unambiguous communication of the meaning of Arithmetical propositions is not possible\footnote{A consequence of Kurt G\"{o}del's formally unchallenged interpretation of his own formal reasoning in his seminal 1931 paper \cite{Go31} on formally undecidable arithmetical propositions.}---are implicit (in the sense that they are neither explicitly accepted nor negated) in the four articles on Formal Proof in the December 2008 Notices of The American Mathematical Society\footnote{Volume 55, Issue 11.}; yet all the authors, curiously, implicitly suggest that effective and unambiguous formalisation of the meaning of Arithmetical propositions, followed by equally effective and unambiguous interpretation of the formalisation in a computational model, \textit{is} possible!

The distinction sought to be emphasised here is between the adequacy\footnote{See, for instance, Wang (\cite{Wa63}, p.5) for a perspective on `adequacy'.}---and limitations---of a formal language when used for expressing abstract concepts precisely, and the adequacy---and limitations---of the language when used to communicate such concepts effectively.

In other words, although a formalisation may be subjectively valid if it can, somehow, be interpreted to mean that which we intended it to express, it is objectively valid only if it communicates that which we intended it to express under any interpretation.

More precisely, a formalisation may express abstract ideas adequately if it is consistent, but it communicates them effectively only if it is categorical.

The thesis underlying the issues addressed in this investigation is that it is the latter---and not merely the former---which needs to be focused upon as an achievable goal in a civilisation that is at the stage of actively developing---and seeking---non-human intelligences.

I shall show how such focus---which should emerge explicitly as the natural motivating factor at some stage in the development of a formal language---appears to be lagely absent so far in the mainstream literature on two of the more basic, formal, first-order languages of mathematics, namely Peano Arithmetic, PA, and Zermelo-Fraenkel set theory, ZF.

I attribute such absence largely to the effect of Nathanson's ``boss" factor in influencing a selective, and subjective, interpretation of Kurt G\"{o}del's formal reasoning in his seminal 1931 paper ``On formally undecidable propositions of Principia Mathematica and related systems I"\footnote{\cite{Go31}.}.

In this paper, G\"{o}del showed the unexpected consequences that follow if we construct a language whose significant elements can be determined effectively (i.e., by mechanical methods only and without use of any ingenuity).

Moreover, G\"{o}del's insistence on using only constructive arguments in the formal results of the paper can, today, be viewed as having prepared the groundwork for Alan Turing's equally famous 1936 paper on computable functions\footnote{\cite{Tu36}.}; the latter, essentially, shows how to translate a first-order version of G\"{o}del's formally defined language of Peano Arithmetic into a machine language (i.e., into configurations that could be effectively recognised, manipulated and transmitted by---and between---mechanical devices).

My thesis---contrary to accepted wisdom---is that, when interpreted constructively, the combination of these two seminal expositions actually yields a means of effective communication between different forms of intelligence that are exposed to a common environment.

\section{Overview and Definitions}

In the first section, `Cohen and the Axiom of Choice', I show that any interpretation of ZF that admits Aristotle's particularisation is not sound.

In the next section, `G\"{o}del and formally undecidable arithmetical propositions', I show that the standard interpretation of PA is not sound, and that PA is \(\omega\)-\textit{in}consistent.

In the section, `Rosser and formally undecidable arithmetical propositions', I show that Rosser's claim in his 1936 paper\footnote{\cite{Ro36}.}---that G\"{o}del's reasoning can be recast to yield a formally undecidable arithmetical proposition without the assumption of \(\omega\)-consistency---is unsustainable, since Rosser implicitly assumes \(\omega\)-consistency in his argument.

In the fourth section, `Turing and a finitary interpretation of PA', I show why an uncritical acceptance of G\"{o}del's interpretation of his own reasoning in his 1931 paper{\cite{Go31}.}---as having established that a formal system of Arithmetic such as PA must have multiple interpretations that are sound, but essentially different---may have led to an equally critical failure to recognise---particularly in view of Turing's 1936 paper\footnote{\cite{Tu36}.}---that PA is categorical (i.e., it essentially admits only one sound interpretation).

In the final section, `Cantor and ZF set theory', I show that although we can relativise PA in ZF so that every PA-theorem relativises as a ZF-theorem, this does not address the issue of whether PA can be consistently extended to ZF (which would be the case if, and only if, ZF has a sound interpretation that is a sound interpretation of PA). I show that this is not the case, and that ZF is inconsistent.

\begin{quote}
\footnotesize{
\textbf{Preliminary Definitions and Comments}

\textbf{Aristotlean particularisation}: This holds that an assertion such as, `There exists an \(x\) such that \(F(x)\) holds'---usually denoted symbolically by `\((\exists x)F(x)\)'---can be validly inferred in the classical logic of predicates from the assertion, `It is not the case that, for any given \(x\), \(F(x)\) does not hold'---usually denoted symbolically by `\(\neg(\forall x)\neg F(x)\)'\footnote{\cite{HA28}, pp.58-59.}.

\textbf{Formal language}: By a ``formal language" I mean a ``formal system" or a ``formal theory" S which meets the following criteria\footnote{Excerpted from \cite{Me64}, p29; however, see also \cite{Sh67}, p.2; \cite{EC89}, p.154.}:

``(1) A countable set of symbols is given as the symbols of S. A finite sequence of symbols of S is called an \textit{expression} of S.

(2) There is a subset of the expressions of S called the set of \textit{well-formed formulas} (abbreviated ``wfs") of S. (There is usually an effective procedure to decide whether a given expression is a wf.)

(3) A set of wfs is set aside and called the set of \textit{axioms} of S. (Most often, one can effectively decide whether a given wf is an axiom, and, in this case, S is called an \textit{axiomatic} theory.)

(4) There is a finite set R\(_{1}\), ..., R\(_{n}\) of relations among wfs, called \textit{rules of inference}. For each R\(_{i}\), there is a unique positive integer \(j\) such that, for every set of \(j\) wfs and each wf \(A\), one can effectively decide whether the given \(j\) wfs are in the relation R\(_{i}\) to \(A\), and, if so, \(A\) is called a \textit{direct consequence} of the given wfs by virtue of one of the rules of inference."

\textbf{Interpretation}: The word ``interpretation" may be used both in its familiar, linguistic, sense, and in a mathematically precise sense; the appropriate meaning is usually obvious from the context. Mathematically, I follow Tarski's definitions of ``interpretation" with respect to the symbolic expressions of a formal language\footnote{\cite{Ta33}.} as excerpted below from Mendelson\footnote{\cite{Me64}, \S2, p49; but see also \cite{Be59}, p.179; \cite{Sh67}, p.61; \cite{BBJ03}; p.102.}:

``An interpretation consists of a non-empty set D, called the domain of the interpretation, and an assignment to each predicate letter \(A_{j}^{n}\) of an \(n\)-place relation in D, to each function letter \(f_{j}^{n}\) of an \(n\)-place operation in D (i.e., a function from D\(_{n}\) into D), and to each individual constant \(a_{i}\) of some fixed element of D. Given such an interpretation, variables are thought of as ranging over the set D, and \(\neg, \rightarrow\), and quantifiers are given their usual meaning. (Remember that an \(n\)-place relation in D can be thought of as a subset of D\(_{n}\), the set of all \(n\)-tuples of elements of D.)"

\textbf{Notation}: I shall henceforth use square brackets to indicate that the contents represent a symbol or a formula of a formal theory, generally assumed to be well-formed unless otherwise indicated by the context.

In other words, expressions inside the square brackets are to be only viewed syntactically as juxtaposition of symbols that are to be formed and manipulated upon strictly in accordance with specific rules for such formation and manipulation---in the manner of a mechanical or electronic device---without any regards to what the symbolism might represent semantically under an interpretation that gives them meaning. 

Moreover, even though the formula `\([R(x)]\)' of a formal Arithmetic may interpret as the arithmetical relation expressed by `\(R^{*}(x)\)', the formula `\([(\exists x)R(x)]\)' need not interpret as the arithmetical proposition denoted by the abbreviation `\((\exists x)R^{*}(x)\).' The latter denotes the phrase `There is some \(x\) such that \(R^{*}(x)\)'. I shall show\footnote{See also \cite{An08f}.} why---as Brouwer had noted\footnote{\cite{Br08}.}---this concept is not always capable of an unambiguous meaning that can be represented in a formal language by the formula `\([(\exists x)R(x)]\)'.

By `expressed' I mean here that the symbolism is simply a short-hand abbreviation for referring to abstract concepts that may, or may not, be capable of a precise `meaning'. Amongst these are symbolic abbreviations which are intended to express the abstract concepts---particularly those of `existence'---involved in propositions that refer to non-terminating processes and infinite aggregates.

\textbf{Provability}: A formula \([F]\) of a formal system S is provable in S (S-provable) if, and only if, there is a finite sequence of S-formulas \([F_{1}], [F_{2}], ..., [F_{n}]\) such that \([F_{n}]\) is \([F]\) and, for all \(1 \leq i \leq n\), \([F_{i}]\) is either an axiom of S or a consequence of the formulas preceding it in the sequence by means of the rules of deduction of S.

\textbf{The first-order Peano Arithmetic (PA)}

\addvspace{+1ex}
\textbf{PA\(_{1}\)}:  \([(x_{1} = x_{2}) \rightarrow ((x_{1} = x_{3}) \rightarrow (x_{2} = x_{3}))]\);

\textbf{PA\(_{2}\)}:  \([(x_{1} = x_{2}) \rightarrow (x_{1}^{\prime} = x_{2}^{\prime})]\);

\textbf{PA\(_{3}\)}:  \([0 \neq x_{1}^{\prime}]\);

\textbf{PA\(_{4}\)}:  \([(x_{1}^{\prime} = x_{2}^{\prime}) \rightarrow (x_{1} = x_{2})]\);

\textbf{PA\(_{5}\)}:  \([( x_{1} + 0) = x_{1}]\);

\textbf{PA\(_{6}\)}:  \([(x_{1} + x_{2}^{\prime}) = (x_{1} + x_{2})^{\prime}]\);

\textbf{PA\(_{7}\)}:  \([( x_{1} \star 0) = 0]\);

\textbf{PA\(_{8}\)}:  \([( x_{1} \star x_{2}^{\prime}) = ((x_{1} \star x_{2}) + x_{1})]\);

\textbf{PA\(_{9}\)}:  For any well-formed formula \([F(x)]\) of PA:

\hspace{+3ex} \([(F(0) \rightarrow (\forall x)(F(x) \rightarrow F(x^{\prime}))) \rightarrow (\forall x)F(x)]\);

\textbf{Modus Ponens in PA}: If [\(A\)] and [\(A \rightarrow B\)] are PA-provable, then so is [\(B\)];

\textbf{Generalisation in PA}: If [\(A\)] is PA-provable, then so is [\((\forall x)A\)].

\textbf{Standard interpretation}: The standard interpretation \(\mathcal{I}_{PA(Standard/Tarski)}\) of PA over the structure \(\mathcal{N}\) is the one in which the logical constants have their `usual' interpretations\footnote{See \cite{Me64}, p.49.} and\footnote{See \cite{Me64}, p.107.}:

(a) the set of non-negative integers is the domain;

(b) the integer 0 is the interpretation of the symbol [0];

(c) the successor operation (addition of 1) is the interpretation of the \([']\) function;

(d) ordinary addition and multiplication are the interpretations of \([+]\) and \([*]\);

(e) the interpretation of the predicate letter \([=]\) is the identity relation.

\textbf{The structure \(\mathcal{N}\)}: The structure of the natural numbers---namely, \{\(N\) (\textit{the set of natural numbers}); \(=\) (\textit{equality}); \('\) (\textit{the successor function}); \(+\) (\textit{the addition function}); \( \ast \) (\textit{the product function}); \(0\) (\textit{the null} \textit{element})\}.

\textbf{Simple consistency}: A formal system S is simply consistent if, and only if, there is no S-formula \([F(x)]\) for which both \([(\forall x)F(x)]\) and \([\neg(\forall x)F(x)]\) are S-provable.

\textbf{\(\omega\)-consistency}: A formal system S is \(\omega\)-consistent if, and only if, there is no S-formula \([F(x)]\) for which, first, \([\neg(\forall x)F(x)]\) is S-provable and, second, \([F(a)]\) is S-provable for any given S-term \([a]\).

\textbf{Soundness}:  A formal system S is sound under an interpretation \(\mathcal{I_{S}}\) if, and only if, every theorem \([T]\) of S translates as `\([T]\) is true under \(\mathcal{I_{S}}\)'.

\textbf{Categoricity}: A formal system S is categorical if, and only if, it has a sound interpretation and any two sound interpretations of S are isomorphic.\footnote{Compare \cite{Me64}, p.91.}

\textbf{Axiom of Choice} (a standard interpretation): Given any set \(S\) of mutually disjoint non-empty sets, there is a set \(C\) containing a single member from each element of \(S\).
}
\end{quote}

\section{Cohen and the Axiom of Choice}

I begin by questioning the belief that a sound interpretation of ZF can appeal to Aristotle's particularisation---a belief that is essential to Paul J.\ Cohen's argument\footnote{In his 1963-64 papers, ``The Independence of the Continuum Hypothesis", \cite{Co63} \& \cite{Co64}.} that the Axiom of Choice is independent of ZF---and show that it admits some curious consequences.

\subsection{Aristotle's particularisation: stronger than the Axiom of Choice}

Now, a fundamental tenet of classical logic---unrestrictedly adopted by formal first-order predicate calculus as axiomatic\footnote{See \cite{Hi25}, p.382; \cite{HA28}, p.48; \cite{Sk28}, p.515; \cite{Be59}, pp.178 \& 218; \cite{Co66}, p.4.}---is Aristotlean particularisation.

This holds that an assertion such as, `There exists an \(x\) such that \(F(x)\) holds'---usually denoted symbolically by `\((\exists x)F(x)\)'---can be validly inferred in the classical logic of predicates from the assertion, `It is not the case that, for any given \(x\), \(F(x)\) does not hold'---usually denoted symbolically by `\(\neg(\forall x)\neg F(x)\)'\footnote{\cite{HA28}, pp.58-59.}.

In his 1908 paper, ``The unreliability of the logical principles"\footnote{\cite{Br08}.}, L.\ E.\ J.\ Brouwer objected to such inference, on the ground that it is invalid as a general logical principle in the absence of a means for constructing some putative object \(a\) that satisfies \(F(a)\).

However, Brouwer's objection did not find favour with the ``bosses" of the day. For instance, David Hilbert dismissively noted in an address on the ``The Foundations of Mathematics" delivered in July 1927 at the Hamburg Mathematical Seminar:

\begin{quote}
\footnotesize{
``\footnote{\cite{Hi27}, p.474.}Brouwer declares [just as Kronecker did in his day] that existence statements, one and all, are meaningless in themselves unless they also contain the construction of the object asserted to exist; for him they are worthless scrip, and their use causes mathematics to degenerate into a game."
}
\end{quote}

A contributory factor for not addressing Brouwer's objection with the seriousness that it merited was, perhaps, the fact that Brouwer formulated it as part of an---unnecessarily stringent---`intuitionist' program that was unable to accommodate the emerging `formalist' perspective of the time.

For instance, Brouwer's rejection of the principle of the excluded middle\footnote{For various informal and formal formulations of this principle see \cite{Br23}, p.335; \cite{Hi25}, p.382; \cite{Hi27}, p.466.} (unwarranted, by the yardstick of Occam's razor) may have merely served to obfuscate the issue.

The reason: Although Brouwer's objection to the validity of Aristotle's particularisation as a general logical principle was justified\footnote{See, for instance, \cite{An08f}.}, it was not for the various reasons cited by him\footnote{See, for instance, \cite{Br23}, pp.336-337.}.

The invalidity is actually due to an ambiguity arising from lack of a specification of the means by which quantified statements---such as `For any given \(x\), \(F(x)\) holds', and `It is not the case that, for any given \(x\), \(F(x)\) holds'---can be validly asserted.

Brouwer failed to recognise that, although the principle of the excluded middle follows from Aristotle's particularisation, the converse is not true.

\subsection{Hilbert's formalisation of Aristotle's particularisation}

Now, in his 1927 address, Hilbert reviewed in detail his axiomatisation of classical Aristotlean predicate logic as a formal first-order \(\varepsilon\)-predicate calculus\footnote{\cite{Hi27}, pp.465-466.}, in which he used a primitive choice-function\footnote{See \cite{Hi25}, p.382.} symbol, `\(\varepsilon\)', for defining the quantifiers `\(\forall\)' and `\(\exists\)' (see Appendix A).

In an earlier address ``On The Infinite", delivered in M\"{u}nster on 4\(^{th}\) June 1925 at a meeting of the Westphalian Mathematical Society, Hilbert had shown that the formalisation proposed by him would adequately express Aristotle's logic of predicates if the \(\varepsilon\)-function was interpreted to yield Aristotlean particularisation\footnote{\cite{Hi25}, pp.382-383; \cite{Hi27}, p.466(1).}.

This, as Hilbert remarked in his 1927 address, was what he had set out to achieve as part of his `proof theory':

\begin{quote}
\footnotesize{
``\footnote{\cite{Hi27}, p.475.}\ldots The fundamental idea of my proof theory is none other than than to describe the activity of our understanding, to make a protocol of the rules according to which our thinking actually proceeds."
}
\end{quote}

What came to be known later as Hilbert's Program\footnote{See, for instance, the Stanford Encyclopedia of Philosophy: \href{http://plato.stanford.edu/entries/hilbert-program/}{Hilbert's Program}.}---which was built upon Hilbert's `proof theory'---can be viewed as, essentially, the subsequent attempt to show that the formalisation was also necessary for communicating Aristotle's logic of predicates effectively and unambiguously under any interpretation of the formalisation.

This goal is implicit in Hilbert's remarks:

\begin{quote}
\footnotesize{
``\footnote{\cite{Hi25}, p.384.}Mathematics in a certain sense develops into a tribunal of arbitration, a supreme court that will decide questions of principle---and on such a concrete basis that universal agreement must be attainable and all assertions can be verified."

``\footnote{\cite{Hi27}, p.475.}\ldots a theory by its very nature is such that we do not need to fall back upon intuition or meaning in the midst of some argument."
}
\end{quote}

\subsection{Aristotle's particularisation and the Axiom of Choice}

The difficulty in attaining this goal constructively along the lines desired by Hilbert---in the sense of the above quotes---lies in the fact that, as Rudolf Carnap emphasised in a 1962 paper, ``On the use of Hilbert's \(\varepsilon\)-operator in scientific theories", the Axiom of Choice is formally derivable as a theorem in a set theory ZF\(_{\varepsilon}\), which is, essentially, Zermelo-Fraenkel set theory where the quantifiers are defined (see Appendix A) in terms of Hilbert's \(\varepsilon\)-function\footnote{\cite{Ca62}, pp.157-158; see also Wang's remarks \cite{Wa63}, pp.320-321.}.

Now, as Thoralf Skolem emphasised in ``Some remarks on axiomatized set theory", delivered in an address before the Fifth Congress of Scandinavian Mathematicians in Helsinki, 4-7 August 1922, the Axiom of Choice is an essentially non-verifiable statement\footnote{\cite{Sk22}, p.300(8).} which does not express any definable content; therefore it cannot be expected to communicate any meaningful information unambiguously under interpretation:

\begin{quote}
\footnotesize{
``So long as we are on purely axiomatic ground there is, of course, nothing special to be remarked concerning the principle of choice (though, as a matter of fact, new sets are \textit{not} generated \textit{univocally} by applications of this axiom); but if many mathematicians---indeed, I believe, most of them---do not want to accept the principle of choice, it is because they do not have an axiomatic conception of set theory at all. They think of sets as given by specification of arbitrary collections; but then they also demand that every set be definable. We can, after all, ask: What does it mean for a set to exist if it can perhaps never be defined? It seems clear that this existence can be only a manner of speaking, which can lead only to purely formal propositions---perhaps made up of very beautiful \textit{words}---about objects \textit{called} sets. But most mathematicians want mathematics to deal, ultimately, with performable computing operations and not to consist of formal propositions about objects called this or that."
}
\end{quote}

Skolem's remarks seem to be directed against Hilbert's explicit focus on viewing the formal consistency of a formal language as adequate for using it as a language of, both, precise expression and effective communication.

Skolem appears to caution that such an approach encourages viewing mathematics as a meaningless ``game" (a caution dismissively referred to by Hilbert in the 1927 address cited above\footnote{\cite{Hi27}, p.474.}) by delinking the development of a formal language, first, from its roots as a formalisation of abstractions that we intend it to express precisely; and, second, from its goal as a means of communicating such abstractions effectively and unambiguously.

As Paul Bernays and Abraham A.\ Fraenkel also caution\footnote{\cite{BF58}.}, a consequence of treating the putative entities postulated by the Axiom of Choice as validly existing in a notional, abstract, space is that:

\begin{quote}
\footnotesize{
``\footnote{\cite{BF58}, p.17, fn.(1).}Lebesgue showed by a concrete example that the distinction between construction and existence may, through the Axiom of Choice, affect elementary (geometrical) problems."
}
\end{quote}

Since Hilbert's formalisation of Aristotle's particularisation can be shown to be stronger than, or equivalent to, the Axiom of Choice of a set theory such as ZF, Aristotle's particularisation too must be viewed as essentially non-verifiable, and nebulous in content, in Skolem's sense!

Brouwer critiqued Hilbert's focus more directly in his 1923 address, delivered on 21\(^{st}\) September 1923, at the annual convention of the Deutsche Mathematiker-Vereinigung in Marburg an der Lahn:

\begin{quote}
\footnotesize{
``\footnote{\cite{Br23}, p.336.}An a priori character was so consistently ascribed to the laws of theoretical logic that these laws, including the principle of the excluded middle, were applied without reservation even in the mathematics of infinite systems and we did not allow ourselves to be disturbed by the consideration that the results obtained in this way are in general not open, either practically or theoretically, to any empirical corroboration. On this basis extensive incorrect theories were constructed, especially in the last half-century. The contradictions that, as a result, one repeatedly encountered gave rise to the \textit{formalistc critique}, a critique which in essence comes to this : the \textit{language accompanying the mathematical mental activity} is subjected to a mathematical examination. To such an examination the laws of theoretical logic present themselves as operators acting on primitive formulas or axioms, and one sets himself the goal of transforming these axioms in such a way that the linguistic effect of the operators mentioned (which are themselves retained unchanged) can no longer be disturbed by the appearance of the linguistic figure of a contradiction. We need by no means despair of reaching this goal,\(^{\tiny{*}}\) but nothing of mathematical value will thus be gained : an incorrect theory, even if it cannot be inhibited by any contradiction that would refute it, is none the less incorrect, just as a criminal policy is none the less criminal even if it cannot be inhibited by any court that would  curb it.

\(^{\tiny{*}}\)For the unjustified application of the principle of excluded middle to properties of well-constructed mathematical systems can never lead to a contradiction \ldots." 
}
\end{quote}

\subsubsection{Any interpretation of ZF which appeals to Aristotle's particularisation is not sound}

The significance of these remarks is seen in the accepted interpretation of Cohen's argument in his 1963-64 papers\footnote{\cite{Co63} \& \cite{Co64}.}.

Cohen's argument is accepted as definitively establishing that the Axiom of Choice is essentially independent of a set theory such as ZF.

Now, this argument---in common with the arguments of many important theorems in standard texts on the foundations of mathematics and logic---appeals to Aristotle's particularisation when interpreting the existential axioms of ZF (or statements about ZF ordinals) in the application of the seemingly paradoxical\footnote{See Skolem's remarks \cite{Sk22}, p295.} (downwards) L\"{o}wenheim-Skolem Theorem\footnote{\cite{Co66}, p.19.} for legitimising putative models of a language (such as the standard model `M'\footnote{\cite{Co66}, p.19 \& p.82.} of ZF and its forced derivative `N'\footnote{\cite{Co66}, p.121.}, in Cohen's argument\footnote{\cite{Co66}, p.83 \& p.112-118.}).

\begin{quote}
\footnotesize{
\textbf{(Downwards) L\"{o}wenheim-Skolem Theorem}\footnote{\cite{Lo15}, p.245, Theorem 6; \cite{Sk22}, p.293.}: If a first-order  proposition is satisfied in any domain at all, then it is already satisfied in a denumerably infinite domain.
}
\end{quote}

Hence Cohen's argument is also applicable to ZF\(_{\varepsilon}\). However, since ZF\(_{\varepsilon}\) proves the Axiom of Choice\footnote{\cite{Ca62}, pp.157-158; see also Wang's remarks \cite{Wa63}, pp.320-321.}, Cohen's argument\footnote{\cite{Co63} \& \cite{Co64}; \cite{Co66}.}---when carried out in ZF\(_{\varepsilon}\)---actually shows that:

\begin{theorem}
: ZF\(_{\varepsilon}\) has no sound model that appeals to Aristotle's particularisation.
\end{theorem}

Now, if ZF has a sound model, then this is also a sound model of ZF\(_{\varepsilon}\). We thus have that:

\begin{corollary}
: ZF has no sound model that appeals to Aristotle's particularisation.
\end{corollary}

We cannot, therefore, conclude that the Axiom of Choice is essentially independent of the axioms of ZF, since none of the models `forced' by Cohen (in his argument for such independence) are sound.

In fact, treating Cohen's argument as a valid `proof' of such independence illustrates the point of Skolem's cautionary remark regarding the significance, and purpose, in developing such formal systems for working mathematicians since, although ZF appeals to Aristotle's particularisation under interpretation, it does not---unlike ZF\(_{\varepsilon}\)---adequately formalise our abstract notions of the logic of predicates as expressed by Aristotle; nor, as Cohen argues, does ZF communicate a categorical interpretation of the Axiom of Choice.

In other words, ZF is neither an adequate language for the precise expression of our subjective abstract concepts \footnote{ZF\(_{\varepsilon}\) being a better candidate.}, nor an adequate language for effectively communicating such concepts objectively.

As Fraenkel and Yehoshua Bar-Hillel remark:

\begin{quote}
\footnotesize{
``\footnote{\cite{FH58}, p.184.}There is no reason to suppose that in a set theory constructed on the basis of an \(\varepsilon\)-calculus the principle of choice would become generally derivable, unless the specific axioms of that set theory contain \(\varepsilon\)-terms themselves."
}
\end{quote}

\subsection{The Continuum Hypothesis}

Moreover, believing in the \textit{essential} independence of the Continuum Hypothesis may be a curiouser consequence of Nathanson's ``boss" factor.

\begin{quote}
\footnotesize{
\textbf{Definition}

\textbf{The Continuum Hypothesis}\footnote{\cite{Co66}, p.67.}: \(c = \aleph_{1}\), where \(c\) is the cardinality of the set of real numbers, and \(\aleph_{1}\)---the immediately larger cardinal\footnote{For a formal definition of the cardinals see \cite{BF58}, p.139; \cite{Me64}, p.184; \cite{Co66}, p.66.} than \(\aleph_{0}\), the cardinality of the natural numbers---is the cardinality of the constructive ordinals.

\textbf{The constructive ordinals}\footnote{See for instance \cite{Be59}, pp.374-375; \cite{Hi25}, p.374.}: The constructive ordinals are \(0_{o}, 1_{o}, 2_{o}, \ldots, \omega, \omega +_{o} 1_{o}, \omega +_{o} 2_{o}, \ldots, \omega.2_{o}, \omega.2_{o} +_{o} 1_{o}, \ldots, \omega^{2_{o}}, \ldots, \omega^{3_{o}}, \ldots, \omega^{\omega}, \ldots, \omega^{\omega^{\omega}}, \ldots <_{o} \epsilon_{0}.\)
}
\end{quote}

For instance, if Ordinal Arithmetic\footnote{See, for instance, \cite{Me64}, p.187.} is consistent, then the members of the set of real numbers between \(0\) and \(1\)---when each real number in the set is expressed in its unique binary form \(0.a_{1}a_{2}\ldots\), with \(a_{i}\) either 0 or 1---are clearly in some (non-algorithmic) 1-1 correspondence with the members of the proper subset of the constructive ordinals, \{\((\omega^{1_{o}}a_{1_{o}} +_{o} \hspace{+.5ex} \omega^{2_{o}}a_{2_{o}} +_{o} \ldots +_{o} \hspace{+.5ex} \omega^{\omega})\)\}.

\begin{quote}
\footnotesize{
(Where, as indicated, each ordinal in the subset is expressed in its unique Cantor normal form; where \(0_{o}, 1_{o}, 2_{o}, \ldots\) denote the finite ordinals that correspond 1-1 to the natural numbers \(0, 1, 2, \ldots\); where \(+_{o}\) and \(<_{o}\) denote ordinal addition and ordinal inequality (corresponding to the relation `less than') respectively; and where \(a_{i_{o}}\) denotes the finite ordinal that corresponds to the natural number denoted by \(a_{i}\).)
}
\end{quote}

Further, the ordinals upto \(\epsilon_{0}\) are constructive, and there is no cardinal number between \(\aleph_{0}\) (which is defined as the cardinality of the set of finite ordinals) and \(\aleph_{1}\) (which is defined as the cardinality of the set of constructive ordinals).

Now, \(c\)---the cardinality of the set of real numbers---is the cardinality of the set of number-theoretic functions of an integral argument whose values are also finite integers\footnote{See \cite{Hi25}, p.384.}; and it is also the cardinality of the set of real numbers between \(0\) and \(1\).

We thus have the curious conclusion that, in every sound interpretation of Ordinal Arithmetic, \(c\) is equal to the interpretation of \(\aleph_{1}\), since \{\((\omega^{1_{o}}a_{1_{o}} +_{o} \hspace{+.5ex} \omega^{2_{o}}a_{2_{o}} +_{o} \ldots +_{o} \hspace{+.5ex} \omega^{\omega})\)\} \(\in \omega^{\omega^{\omega}}\) and \(\omega^{\omega^{\omega}} <_{o} \epsilon_{0}\)!

Now, if Ordinal Arithmetic is consistent then, by G\"{o}del's completeness theorem, this implies that the Continuum Hypothesis, too, is expressible as a theorem in it\footnote{This, essentially, was Hilbert's argument and thesis in \cite{Hi25}, p.384.}!

\begin{quote}
\footnotesize{
\textbf{G\"{o}del's Completeness Theorem}: In any first-order predicate calculus, the theorems are precisely the logically valid well-formed formulas (i.\ e.\ those that are true in every model of the calculus).}
\end{quote}

This, of course, is at variance with the accepted belief that the Continuum Hypothesis is not provable in Ordinal Arithmetic; a belief that---in no small measure---can be attributed to Nathanson's ``boss" factor endorsing the point of view Cohen expressed at the conclusion of his lectures on ``Set Theory and the Continuum Hypothesis", delivered at Harvard University in the spring term of 1965. 

Thus, whilst, essentially, considering the argument sketched out above, Cohen remarked:

\begin{quote}
\footnotesize{
``\footnote{\cite{Co66}, p.151.}We close with the observation that the problem of CH is not one which can be avoided by not going up in type to sets of real numbers. A similar undecidable problem can be stated using only the real numbers. Namely, consider the statement that every real number is constructible by a countable ordinal. Instead of speaking of countable ordinals we can speak of suitable subsets of \(\omega\). The construction \(\alpha \rightarrow F_{\alpha}\) for \( \alpha \leq \alpha_{0}\), where \(\alpha_{0}\) is countable, can be completely described if one merely gives all pairs \((\alpha, \beta)\) such that \(F_{\alpha} \in F_{\beta}\). This in turn can be coded as a real number if one enumerates the ordinals. In this way one only speaks about real numbers and yet has an undecidable statement in ZF. One cannot push this farther and express any of the set-theoretic questions that we have treated as statements about integers alone. Indeed one can postulate as a rather vague article of faith that any statement in arithmetic is decidable in ``normal" set theory, i.e., by some recognizable axiom of infinity. This is of course the case with the undecidable statements of G\"{o}del's theorem which are immediately decidable in higher systems."
}
\end{quote}

Curiously, Cohen appears to assert here that if ZF is consistent, then we can `see' that the Continuum Hypothesis is subjectively true for the integers under some model of ZF, but---along with the Generalised Continuum Hypothesis---we cannot objectively `assert' it to be true for the integers since it is not provable in ZF, and hence not true in all models of ZF.

However, by this argument, G\"{o}del's undecidable arithmetical propositions, too, can be `seen' to be subjectively true for the integers in the standard model of PA, but cannot be `asserted' to be true for the integers since the statements are not provable in an \(\omega\)-consistent PA, and hence they are not true in all models of an \(\omega\)-consistent PA!

\begin{quote}
\footnotesize{
The latter is, essentially, John Lucas' well-known G\"{o}delian argument\footnote{\cite{Lu61}.}, forcefully argued by Roger Penrose in his popular expositions, `Shadows of the Mind'\footnote{\cite{Pe94}.} and `The Emperor's New Mind'\footnote{\cite{Pe90}.}.

As I have argued in The Reasoner\footnote{\cite{An07a}; \cite{An07b}; \cite{An07c}.}, the argument is plausible, but unsound.

Moreover, the argument needs to be viewed more appropriately as an instance of Nathanson's ``boss" factor, since it is based on a misinterpretation---of what G\"{o}del actually proved formally in his 1931 paper---for which neither Lucas nor Penrose ought to be taken to account\footnote{\cite{An07b}; \cite{An07c}.}.
}
\end{quote}

The distinction sought to be drawn by Cohen is curious, since we have shown that his argument---which assumes that sound interpretations of ZF can appeal to Aristotle's particularisation---actually establishes that sound interpretations of ZF cannot appeal to Aristotle's particularisation; just as I show in the next section that G\"{o}del's argument\footnote{In \cite{Go31}, p.24, Theorem VI.} actually establishes that sound interpretations of PA, too, cannot appeal to Aristotle's particularisation.

Loosely speaking, the cause of the undecidability of the Continuum Hypothsis, and of the Axiom of Choice, in ZF as shown by Cohen, and that of G\"{o}del's undecidable proposition in Peano Arithmetic, is common; it is interpretation of the existential quantifier under an interpretation as Aristotlean particularisation.

In Cohen's case, such interpretation is made explicitly and unrestrictedly in the underlying predicate logic\footnote{\cite{Co66}, p.4.} for ZF; in G\"{o}del's case it is made explicitly---but in a restricted sense to avoid attracting intuitionistic objections---through his specification of the weaker property of \(\omega\)-consistency\footnote{\cite{Go31}, pp.23-24.} (which is a consequence of Aristotle's particularisation) for his formal system P of Peano Arithmetic.

\section{G\"{o}del and formally undecidable arithmetical propositions}

In this section I show that a closer scrutiny of G\"{o}del's formal reasoning in his 1931 paper{\cite{Go31}.} actually implies that PA is \(\omega\)-\textit{in}consistent; and that, as a consequence, the standard interpretation \(\mathcal{I}_{PA(Standard/Tarski)}\) of PA---which, too, appeals to Aristotle's particularisation---is also not sound.

\subsection{The significance of Hilbert's formalisation for Tarski's standard definitions of the satisfaction, and truth, of the formulas of a formal language, under an interpretation}

Now, the larger significance of Hilbert's formalisation of Aristotle's particularisation is that---in formal languages that prefer the more familiar `\([\forall]\)' as a primitive symbol to Hilbert's more logical choice function `\([\varepsilon]\)'---it is Hilbert's formalisation of Aristotle's particularisation\footnote{\cite{Hi27}, pp.465-466.} that implicitly gives formal legitimacy to Alfred Tarski's standard definitions of the satisfaction, and truth, of the formulas of a formal language under an interpretation\footnote{\cite{Ta33}; see also \cite{Me64}, pp.49-53.}, since these definitions faithfully mirror the particular interpretation of Hilbert's formalisation that appeals to Aristotle's particularisation.

The reason: Under Tarski's definitions, the formally defined logical constant `\([\exists]\)' in an occurence such as `\([(\exists x) \ldots ]\)'---which is formally defined in terms of the primitive (undefined) logical constant `\([\forall]\)' as `\([\neg(\forall x)\neg \ldots ]\)'---\textit{always} appeals to an interpretation such as `There is some \(x\) such that \ldots' in \textit{any} formal first-order mathematical language\footnote{See, for instance, \cite{Me64}, p.52(ii).}.

In other words, Tarski's definitions ensure that, if the first-order predicate calculus of a first-order mathematical language admits quantification, then any putative model of the language \textit{must} interpret existential quantification as Aristotle's particularisation\footnote{See, for example, Cohen's standard model `M' (\cite{Co66}, p.19 \& p.82) of ZF, and its forced derivative `N' (\cite{Co66}, p.121), whose existence is legitimised by appeal to the (downward) L\"{o}wenheim-Skolem theorem (\cite{Lo15}, p.245, Theorem 6; \cite{Sk22}, p.293) which, in turn, appeals to Aristotle's particularisation (\cite{Co66}, p.4).}.

\subsection{Consequences of favouring Tarski's interpretation: the standard interpretation of the \\ Peano Arithmetic PA is not sound}

I now consider some consequences of selecting such a strong interpretation---i.e., one which favours Aristotle's particularisation---for the standard interpretation \(\mathcal{I}_{PA(Standard/Tarski)}\) of the formal Peano Arithmetic PA over the structure \(\mathcal{N}\), since \(\mathcal{I}_{PA(Standard/Tarski)}\) appeals to Tarski's definitions.

For instance, if we accept the logical validity of such interpretation, then \(\mathcal{I}_{PA(Standard/Tarski)}\) is sound (i.e., every PA-theorem interprets as true under \(\mathcal{I}_{PA(Standard/Tarski)}\)).

Further, if \(\mathcal{I}_{PA(Standard/Tarski)}\) is sound, then PA is \(\omega\)-consistent (i.e., we cannot have a PA-formula \([F(x)]\) such that \([F(n)]\) is PA-provable for any given PA-numeral \([n]\), and \([\neg (\forall x)F(x)]\) is also PA-provable).

\subsection{The significance of \(\omega\)-consistency: Hilbert's program}

The significance of \(\omega\)-consistency can be traced back to Hilbert's program. Thus, as part of his program for giving mathematical reasoning a finitary foundation, Hilbert\footnote{cf.\ \cite{Hi30}, pp.485-494.} proposed an \(\omega\)-rule as a finitary means of extending a Peano Arithmetic, such as P, to a possible completion (i.e. to logically showing that, given any arithmetical proposition, either the proposition, or its negation, is formally provable from the axioms and rules of inference of the extended Arithmetic).

\begin{quote}
\textbf{Hilbert's \(\omega\)-Rule}: If it is proved that the P-formula [\(F(x)\)] interprets as a true numerical formula for each given P-numeral [\(x\)], then the P-formula [\((\forall x)F(x)\)] may be admitted as an initial formula (\textit{axiom}) in P.
\end{quote}

Now, if we meta-assume Hilbert's \(\omega\)-rule for P, then it follows that, if P is consistent, then there is no P-formula [\(F(x)\)] for which, first, [\(\neg(\forall x)F(x)\)] is P-provable and, second, [\(F(n)\)] is P-provable for any given P-numeral [\(n\)]. In other words, if we meta-assume Hilbert's \(\omega\)-rule for P, then a consistent P is necessarily \(\omega \)-consistent.

\subsection{Aristotle's particularisation implies that PA has non-standard models}

Now, in his seminal 1931 paper\footnote{\cite{Go31}.}, G\"{o}del showed that if a Peano Arithmetic such as his formal system P is \(\omega\)-consistent, then it is incomplete\footnote{\cite{Go31}, Theorem VI, p.24.} (in the sense that he could constructively define a P-formula \([R(x)]\) such that neither \([(\forall x)R(x)]\) nor \([\neg(\forall x)R(x)]\) are P-provable\footnote{\cite{Go31}, p.25(1) \& p.26(2).}).

G\"{o}del concluded that an \(\omega\)-consistent P must, therefore, have a consistent, but \(\omega\)-\textit{in}consistent, extension P\('\), obtained by adding \([\neg(\forall x)R(x)]\) as an axiom to P\footnote{\cite{Go31}, p.27.}.

Since G\"{o}del's argument holds for PA, we thus have that:

\begin{quote}
(a) If Hilbert's formalisation of Aristotlean predicate logic has a sound interpretation when the \(\varepsilon\)-function is interpreted to appeal to Aristotlean particularisation, then the standard interpretation \(\mathcal{I}_{PA(Standard/Tarski)}\) of PA is sound;

(b) If \(\mathcal{I}_{PA(Standard/Tarski)}\) is sound then PA is \(\omega\)-consistent;

(c) If PA is \(\omega\)-consistent then it has a non-standard model containing non-natural numbers that satisfy PA (under Tarski's definitions of the satisfiability and truth of the propositions of a formal language under an interpretation), but which are not definable in it!
\end{quote}

\subsection{A consistent PA is \textit{not} \(\omega\)-consistent}

However, G\"{o}del also showed that if P is consistent and \([(\forall x)R(x)]\) is assumed P-provable, then \([\neg(\forall x)R(x)]\) is P-provable\footnote{This follows from G\"{o}del's argument in \cite{Go31}, p.26(1).}. By G\"{o}del's definition of P-provability, it follows that there is a finite sequence \([F_{1}], \ldots, [F_{n}]\) of P-formulas such that \([F_{1}]\) is \([(\forall x)R(x)]\), \([F_{n}]\) is \([\neg(\forall x)R(x)]\), and, for \(2\leq i \leq n\), \([F_{i}]\) is either a P-axiom or a logical consequence of the preceding formulas in the sequence by the rules of inference of P.

Now, a proof sequence of P necessarily interprets as a sound deduction sequence under any sound interpretation of P. It follows that we cannot have a sound interpretation of P under which \([(\forall x)R(x)]\) interprets as true and \([\neg(\forall x)R(x)]\) as false.

Since both \([(\forall x)R(x)]\) and \([(\forall x)R(x)]\) are closed P-formulas, it follows that the P-formula \([(\forall x)R(x) \rightarrow \neg(\forall x)R(x)]\) interprets as true under every sound interpretation of P. By G\"{o}del's completeness theorem \([(\forall x)R(x) \rightarrow \neg(\forall x)R(x)]\) is, therefore, P-provable; whence \([\neg(\forall x)R(x)]\) is P-provable.

Since G\"{o}del also showed that, if P is consistent, then \([R(n)]\) is P-provable for any given P-numeral \([n]\)\footnote{\cite{Go31}, p.26(2).}, it follows that P is \textit{not} \(\omega\)-consistent (The detailed argument is given in Appendix B).

Since G\"{o}del's argument holds in PA, we further have that:

\begin{theorem}
: A consistent PA is \textit{not} \(\omega\)-consistent.
\end{theorem}

Thus G\"{o}del's Theorem VI of his 1931 paper\footnote{\cite{Go31}, pp.24-26.} is vacuously true, and does not establish the existence of a formally undecidable proposition in P.

(The above argument also applies to J.\ Barkley Rosser's extension of G\"{o}del's reasoning\footnote{\cite{Ro36}.}, since---as I show in the next section---Rosser's argument implicitly appeals to Aristotle's particularisation; thus, despite his claim of having assumed only simple consistency for P, Rosser's argument also presumes---albeit implicitly---that P is \(\omega\)-consistent.)

\subsubsection{The standard interpretation \(\mathcal{I}_{PA(Standard/Tarski)}\) of PA is not sound}

It also follows that a sound interpretation of PA, too, cannot appeal to Aristotle's particularisation, and so:

\begin{theorem}
: The interpretation \(\mathcal{I}_{PA(Standard/Tarski)}\) of PA is not sound.
\end{theorem}

In other words---as Brouwer had noted\footnote{\cite{Br08}.}---the phrase, `There is some \(x\) such that \(R^{*}(x)\)', is not always capable of an unambiguous meaning that can be represented in a formal language by the formula `\([(\exists x)R(x)]\)'.

\subsubsection{PA is categorical}

Nevertheless, as I show in a subsequent section, Turing's seminal 1936 paper on computable numbers\footnote{\cite{Tu36}.} actually admits a sound, finitary, interpretation of an \(\omega\)-\textit{in}consistent PA that is categorical. 

\subsection{G\"{o}del's interpretation of his own formal reasoning: Another case of Nathanson's ``boss" factor}

Meanwhile, I briefly review some factors that appear to have influenced the acceptance of the standard model \(\mathcal{I}_{PA(Standard/Tarski)}\) of PA as sound.

Now, historically it was in the course of an investigation\footnote{Apparently motivated by Hilbert's Program.} intended to capture the essential concepts of Arithmetic---sought to be expressed by the Peano Postulates---in a formal language P that G\"{o}del discovered\footnote{\cite{Go31}, pp.5-6.} he could construct a proposition \([(\forall x)R(x)]\) in P such that, if P is assumed \(\omega\)-consistent, then neither \([(\forall x)R(x)]\) nor \([\neg(\forall x)R(x)]\) could be proven from the axioms of P by the rules of inference of P\footnote{See \cite{Go31}, p.24, Theorem VI.}!

Moreover, G\"{o}del constructed his `undecidable' P-proposition \([(\forall x)R(x)]\) in such a way that its `standard' interpretation was an assertion about the natural numbers that was undeniably true (since it could constructively be shown to be a logical consequence of the Peano Postulates).

\begin{quote}
\footnotesize{
``\footnote{\cite{Go31}, p.26.}One can easily convince oneself that the proof we have just given is constructive \ldots For, all the existential assertions occuring in the proof rest upon Theorem V, which, as is easy to see, is intuitionistically unobjectionable". 
}
\end{quote}

\subsection{G\"{o}del's interpretation of his own formal reasoning implies that perfect communication is not theoretically possible}

G\"{o}del showed, further, that this state of affairs would persist for any \(\omega\)-consistent language that seeks to formally express all our true propositions about the natural numbers that are sought to be captured by the Peano Postulates using classical first-order logic.

Now, as noted earlier, G\"{o}del interpreted his reasoning as showing not only that there is more than one way of interpreting any \(\omega\)-consistent language which seeks to formally express all our true propositions about the natural numbers, but that, in some of these `non-standard' interpretations, our formally undecidable but intuitively `true' (i.e., true under the standard interpretation of PA) proposition would be seen as false!

\begin{quote}
\footnotesize{
``\footnote{\cite{Go31}, p.27.}If one adjoins \(Neg(17Gen \hspace{+.5ex} r)\) to \(\kappa\), then one obtains a consistent, but not an \(\omega\)-consistent, class of FORMULAS \(\kappa'\)."
}
\end{quote}

So, even though G\"{o}del's P-unprovable formula \([(\forall x)R(x)]\) would interpret as true under the `standard' interpretation of a consistent P, we could add the P-unprovable formula \([\neg(\forall x)R(x)]\) to P and obtain a consistent language P\('\) such that, under a sound interpretation of P\('\) (which G\"{o}del implicitly implies as existing), the P-formula \([(\forall x)R(x)]\) would interpret as `false'!

In other words---according to G\"{o}del's interpretation of his own formal reasoning in his 1931 paper\footnote{\cite{Go31}.}---even if we are able to represent all our true propositions about the natural numbers in a formal Arithmetic \(A\), there will always be some of our true propositions that will not be known to be true to someone who has access only to the formal expressions of these propositions in \(A\), and who tries to `decode' their meaning in order to determine whether the propositions can be held to be objectively true or not.

This implies that perfect communication is not theoretically possible in such languages.

Now, first-order Peano Arithmetic PA (which is accepted as the most constructive and faithful formalisation of the Peano Postulates) is the foundation for the most fundamental and unequivocal of our machine languages---namely those that enable mechanical and electronic artefacts to `talk' to each other without any perceivable ambiguity.

Thus, acceptance of G\"{o}del's interpretation of his own reasoning  as definitive by the ``bosses" of the day entails acceptance of the conclusion that there are theoretical limits on unambiguous and effective communication.

However, this acceptance is based on a selective interpretation of G\"{o}del's remarks at the conclusion of his Theorem VI\footnote{\cite{Go31}, p.28.}; remarks which appear curiously ambiguous both on the significance of the assumption of \(\omega\)-consistency for his system P of Arithmetic, and on the formal systems in which his arguments are valid.

\begin{quote}
\footnotesize{
``\footnote{\cite{Go31}, p.28.}In the proof of Theorem VI no properties of the system P were used other than the following:

1. The class of axioms and the rules of inference (i.e. the relation ``immediate consequence"') are recursively definable (when the primitive symbols are replaced in some manner by natural numbers).

2. Every recursive relation is definable within the system P (in the sense of Theorem V).

Hence, in every formal system which satisfies assumptions 1, 2 and is \(\omega\)-consistent, there exist undecidable propositions of the form \((x)F(x)\), where \(F\) is a recursively defined property of the natural numbers, and likewise in every extension of such a system by a recursively definable \(\omega\)-consistent class of axioms. To the systems which satisfy assumptions 1, 2 belong, as one can easily confirm, the Zermelo-Fraenkel and the v. Neumann axiom systems for set theory, and, in addition, the axiom system for number theory which consists of Peano's axioms, recursive definitions (according to schema (2)) and the logical rules. \ldots"

``\footnote{\cite{Go31}, p.28, footnote 48a.}The true reason for the incompleteness which attaches to all formal systems of mathematics lies, as will be shown in Part II of this paper, in the fact that the formation of higher and higher types can be continued into the transfinite (c.f., D. Hilbert, `\"{U}ber das Unendliche', Math. Ann. 95, p. 184), while, in every formal system, only countable many are available. Namely, one can show that the undecidable sentences which have been constructed here always become decidable through adjunction of suitable high types (e.g. of the type \(\omega \) to the system \(P\). A similar result also holds for the axiom systems of set theory."
}
\end{quote}

Thus, G\"{o}del's remarks implicitly suggest that all the formal systems of Peano Arithmetic cited by him can be assumed \(\omega\)-consistent.\footnote{\cite{Go31}, p.28. That this belief is both current and widespread is evidenced in the following entry (as of 27th July 2008) under the  \href{http://plato.stanford.edu/entries/self-reference/}{Self-Reference section of the on-line Stanford Encyclopedia of Philosophy (\S2.4)}, where it is remarked that: ``G\"{o}del's theorem can be interpreted as demonstrating a limitation in what can be achieved by purely formal procedures. It says that if first-order arithmetic is \(\omega\)-consistent (\textit{which it is believed to be}), then there must be arithmetical sentences that can neither be proved nor disproved by the formal procedures of first-order arithmetic".}

However, as I have shown above, first-order Peano Arthmetic PA is actually \(\omega\)-\textit{in}consistent!

Moreover, once we admit such interpretations, we must be prepared to face the disquieting---though, prima facie, spiritually implausible---spectre of the possibility of idealogical conflicts with differing intelligences---and not merely differing faiths---over why (as at present), or even whether, a particular interpretation can lay claim to being the preferred `standard' interpretation of a formal system such as PA!

\section{Rosser and formally undecidable arithmetical propositions}

Of course, since every \(\omega\)-consistent system is necessarily simply consistent, G\"{o}del's reasoning per se is significant only if there is an \(\omega\)-consistent language that seeks to formally express all our true propositions about the natural numbers.

Now, the issue of whether there is an \(\omega\)-consistent system of Arithmetic at all appears to have been treated as inconsequential\footnote{See, for instance, \cite{Be59}, p.595; \cite{Wa63}, p.19 (Theorem 3) \& p.25; \cite{Me64}, p.144; \cite{Sh67}, p.132 (Incompleteness Theorem); \cite{EC89}, p.215; \cite{BBJ03}, p.224 (G\"{o}del's first incompleteness theorem).} following Rosser's 1936 paper\footnote{\cite{Ro36}.}, in which he claimed that G\"{o}del's reasoning can be recast to arrive at his intended result without the assumption of \(\omega\)-consistency.

Although both G\"{o}del's proof and Rosser's argument are complex, and not easy to unravel, the former has been extensively analysed and its various steps formally validated\footnote{Possibly because G\"{o}del's remarkably self-contained 1931 paper---it neither contained, nor needed, any formal citations---remains unsurpassed in mathematical literature for thoroughness, clarity, transparency and soundness of exposition.} in a number of expositions of G\"{o}del's number-theoretic reasoning\footnote{For instance \cite{Me64}, p.143; \cite{EC89}, p.210-211.}.

\subsection{Rosser's argument: A further case of Nathanson's ``boss" factor}

In sharp contrast, Rosser's widely cited\footnote{\cite{Be59}, pp.393-395 (which focuses on Rosser's argument, and treats G\"{o}del's proof of his Theorem VI (\cite{Go31}, p.24) as a, secondary, weaker result); \cite{Wa63}, p.337; \cite{Me64}, pp.144-147; \cite{Sh67}, p.232 (interestingly, this introductory text contains \textit{no} reference to G\"{o}del or to his 1931 paper!); \cite{EC89}, p.215; \cite{Sm92}, p.81; \cite{BBJ03}, p.226 (this introductory text, too, focuses on Rosser's argument, and treats G\"{o}del's argument as more of a historical curiosity!).} argument does not appear to have received the same critical scrutiny, and---lending force to Nathanson's observation---its number-theoretic expositions generally remain either implicit or sketchy\footnote{See \cite{Be59}, p.593-595; \cite{Wa63}, p.337; \cite{Sh67}, p.232; \cite{Rg87}, p.98; \cite{EC89}, p.217, Ex.2; \cite{Sm92}, p.81; \cite{BBJ03}, p.226.}.

\subsubsection{The significance of Rosser's claim}

The significance of Rosser's argument is that its acceptance as a `proof' lends legitimacy to arguments which implicitly assume that the existence of undecidable propositions in PA must follow---without addressing the significance or consequences of \(\omega\)-consistency---from the establishment of undecidable propositions in first-order systems \textit{other} than PA that are capable of formalising the Peano Postulates.\footnote{For instance \cite{Sh67} on p.232, and \cite{BBJ03} on pp.224 and 226.}

If Rosser's argument had not been accepted so readily---and so uncritically---as capable of conversion to a formal proof by the ``bosses" of the day\footnote{Including G\"{o}del!}, surely more attention would have been given to the significance of G\"{o}del's explicit assumption of \(\omega\)-consistency for his system of Peano Arithmetic; to its relation to Hilbert's \(\omega\)-rule; to its deeper relation to the Hilbert-Brouwer dispute over the interpretation of the existential quantifier of a formal language that appeals to Aristotle's particularisation; and to the possibility of its falsity.

However, acceptance of Rosser's argument as valid blurred the critical distinction between G\"{o}del's original, constructive\footnote{See \cite{Go31}, p.26.}, number-theoretic, argument---which appeals only to the consistency of first-order formalisations of Peano Arithmetic such as PA---and arguments that:

(i) seek to validate G\"{o}del's Theorem VI of his 1931 paper (such as Rosser's claimed extension of it) by reasoning that appeals to arguments\footnote{Such as \cite{Be59}, p.594; \cite{Wa63}, p.337; \cite{EC89}, p.206 \& p.215; \cite{Sm92}, p.36; \cite{BBJ03}, p.224.} which presume upon the soundness of the `standard' interpretation \(\mathcal{I}_{PA(Standard/Tarski)}\) of PA in an attempt to conclude these results as corollaries of arguments that by-pass the constructive number-theoretic constraints within which G\"{o}del's arguments were originally derived;

(ii) seek to validate G\"{o}del's Theorem VI of his 1931 paper (and Rosser's claimed extension of it) by reasoning that appeals to the consistency of formalisations of Cantor's set theory---such as ZF---in an attempt to conclude these results as corollaries of more general, non-constructive, set-theoretical arguments that by-pass the constructive number-theoretic constraints within which G\"{o}del's arguments were originally derived\footnote{For instance \cite{Sh67}, p.232; \cite{EC89}, p.215.}.

Striking examples of (i) are Hao Wang's and E.\ W.\ Beth's independent, sketchy, formalisations of Rosser's informal argument.

\subsubsection{Wang's outline of Rosser's argument}

Wang, for instance, states that\footnote{\cite{Wa63}, p.337.} from the formal provability of:

\vspace{+1ex}
(i) \(\neg(x)(B(x, \overline{q}) \supset (Ey)(y \leq x \hspace{+.5ex} \& \hspace{+.5ex} B(y, n(\overline{q}))))\)

\vspace{+1ex}
\noindent in his formal system of first-order Peano Arithmetic Z, we may infer the formal provability of:

\vspace{+1ex}
(ii) \((Ex)(B(x, \overline{q}) \hspace{+.5ex} \& \hspace{+.5ex} \neg(Ey)(y \leq x \hspace{+.5ex} \& \hspace{+.5ex} B(y, n(\overline{q}))))\)

\vspace{+1ex}
However, the inference (ii) from (i) only follows if we assume that the following deduction is valid for some \(\overline{j}\):

\begin{quote}
\(\neg(x)(B(x, \overline{q}) \supset (Ey)(y \leq x \& B(y, n(\overline{q}))))\)

\(B(\overline{j}, \overline{q}) \hspace{+.5ex} \& \hspace{+.5ex} \neg(Ey)(y \leq \overline{j} \hspace{+.5ex} \& \hspace{+.5ex} B(y, n(\overline{q})))\)

\((Ex)(B(x, \overline{q}) \hspace{+.5ex} \& \hspace{+.5ex} \neg(Ey)(y \leq x \hspace{+.5ex} \& \hspace{+.5ex} B(y, n(\overline{q}))))\)
\end{quote}

Now, this deduction assumes that the `standard' interpretation of PA, when applied to Wang's Peano Arithmetic Z, is sound---which is not the case!

\begin{quote}
\footnotesize{
Although Wang does not explicitly define the interpretation of the formal Z-formula `\((Ex)F(x)\)' as `There is some \(x\) such that \(F(x)\)', this interpretation appears implicit in his discussion and definition of `\((Ev)A(v)\)' in terms of Hilbert's \(\varepsilon\)-function\footnote{\cite{Wa63}, p.315(2.31); see also p.10 \& pp.443-445.} as a property of the underlying logic of Wang's Peano Arithmetic Z, and is obvious in the above argument.

In other words Wang implicitly implies that the interpretation of existential quantification cannot be specific to any particular interpretation of a formal mathematical language, but must necessarily be determined by the predicate calculus that is to be applied uniformly to all the mathematical languages in question.
}
\end{quote}

\subsubsection{Beth's outline of Rosser's argument}

Similarly, in his outline of a formalisation of Rosser's argument, Beth implicitly concludes\footnote{\cite{Be59}, p.594 (ij).} that from the formal provability of:

\vspace{+1ex}
(i) \(\neg(q)[G_{1}(m^{0}, q, m^{0}) \rightarrow (s)\{B(s, q) \rightarrow (Et)[t \leq s \hspace{+.5ex} \& \hspace{+.5ex} (Er)\{H(q, r) \hspace{+.5ex} \& \hspace{+.5ex} B(t, r)\}]\}]\)

\vspace{+1ex}
\noindent in his formal system of first-order Peano Arithmetic P, we may infer the formal provability of:

\vspace{+1ex}
(ii) \((Eq)[G_{1}(m^{0}, q, m^{0}) \hspace{+.5ex} \& \hspace{+.5ex} (s)\{B(s, q) \hspace{+.5ex} \& \hspace{+.5ex} (t)[t \leq s \rightarrow (r)\{H(q, r) \rightarrow \overline{B(t, r)}\}]\}]\)

\vspace{+1ex}
However, the inference (ii) from (i) only follows if we assume that the following deduction is valid for some \(\overline{j}\):

\begin{quote}
\(\neg(q)[G_{1}(m^{0}, q, m^{0}) \rightarrow (s)\{B(s, q) \rightarrow (Et)[t \leq s \hspace{+.5ex} \& \hspace{+.5ex} (Er)\{H(q, r) \hspace{+.5ex} \& \hspace{+.5ex} B(t, r)\}]\}]\)

\(G_{1}(m^{0}, \overline{j}, m^{0}) \hspace{+.5ex} \& \hspace{+.5ex} (s)\{B(s, \overline{j}) \hspace{+.5ex} \& \hspace{+.5ex} (t)[t \leq s \rightarrow (r)\{H(\overline{j}, r) \rightarrow \overline{B(t, r)}\}]\}\)

\((Eq)[G_{1}(m^{0}, q, m^{0}) \hspace{+.5ex} \& \hspace{+.5ex} (s)\{B(s, q) \hspace{+.5ex} \& \hspace{+.5ex} (t)[t \leq s \rightarrow (r)\{H(q, r) \rightarrow \overline{B(t, r)}\}]\}]\)
\end{quote}

\vspace{+1ex}
This deduction, again, assumes that the `standard' interpretation of PA, when applied to Beth's Peano Arithmetic P, is sound---which is not the case!

\begin{quote}
\footnotesize{
In this case, Beth explicitly defines the interpretation of the formal P-formula `\((Ex)\)' as `There is a value of \(x\) such that'\footnote{\cite{Be59}, p.178.}.

Thus Beth, too, implies that the interpretation of existential quantification in formalised axiomatics cannot be specific to any particular interpretation of a formal mathematical language, but must necessarily be determined by the predicate calculus that is to be applied uniformly to all the mathematical languages in question.
}
\end{quote}

\subsection{Rosser's argument does not support his claim}

Now, Rosser's claim in his `extension'\footnote{\cite{Ro36}.} of G\"{o}del's argument\footnote{\cite{Go31}.} is that, whereas G\"{o}del's argument assumes that his Peano Arithmetic, P, is \(\omega\)-consistent, Rosser's assumes only simple consistency.

P is defined as \(\omega\)-consistent if, and only if, there is no P-formula \( [R(x)] \) for which \( [\neg(\forall x)R(x)] \) is P-provable whilst \( [R(n)] \) is P-provable for any given numeral \( [n] \) of P.

However, Rosser's original argument (sketch) also appears to presume---albeit implicitly---that the system of Peano Arithmetic in question is \(\omega\)-consistent.

For instance, the concluding deduction in Rosser's reasoning \textit{implicitly} presumes that if, for any given natural number \(n\), the formula in G\"{o}del's Peano Arithmetic P whose G\"{o}del-number is:

\(Neg(Sb(r \begin{array}{c} u \\ Z(n) \end{array} \begin{array}{c} v \\ Z(a) \end{array} ))\)

\noindent is P\(_{\kappa}\)-provable\footnote{Notation (due to G\"{o}del): By `P\(_{\kappa}\)-provable' we mean provable from the axioms of PA and an arbitrary class of PA-formulas \(\kappa\)---including the case where \(\kappa\) is empty---by the rules of deduction of PA.} under the given premises, we may conclude that, if P is simply consistent, then the P-formula whose G\"{o}del-number is:

\(uGen(Neg(Sb(r \begin{array}{c} v \\ Z(a) \end{array} )))\)

\noindent is also P\(_{\kappa}\)-provable.

Rosser essentially seems to reason here that since the P-formula \([\neg R(n, a)]\)---where \([a]\) is a specific P-numeral---is P\(_{\kappa}\)-provable for any given P-numeral \([n]\), we may conclude that the P-formula \([(\forall u)\neg R(u, a)]\) is P\(_{\kappa}\)-provable. This would presume, however, that P is \(\omega\)-consistent.

\subsection{Where Rosser's argument \textit{presumes} \(\omega\)-consistency}

I consider a more formal expression\footnote{eg.\ \cite{Me64}, p.145, Proposition 3.32.} of Rosser's argument in Appendix C, and highlight where it implicitly presumes that P is \(\omega\)-consistent.

\section{Turing and a finitary interpretation of PA}

In this section I show why an uncritical acceptance of G\"{o}del's interpretation of his own reasoning in\footnote{\cite{Go31}.}---as having established that a formal system of Arithmetic such as PA must have multiple interpretations that are sound, but essentially different---may have led to a critical failure to recognise (particularly in view of Turing's 1936 paper\footnote{\cite{Tu36}.}) that PA is categorical (i.e., that PA essentially admits only one sound interpretation).

The significance of such recognition is that it would explain why we unhesitatingly entrust our lives each moment of each day to mechanical and electronic artefacts whose reliability is essentialy founded on the ability of PA to admit unambiguous and effective communication. It would also imply that we can, in principle, communicate perfectly with technologically advanced extra-terrestrial intelligences.

\subsection{Turing-computability: A finitary interpretation of PA}

I consider a weakened, finitary, interpretation \(\mathcal{I}_{PA(\beta rouwer/Turing)}\) of an \(\omega\)-\textit{in}consistent PA which avoids appealing to Aristotlean particularisation in the interpretation of the existential quantifier, and which is actually implicit in Turing's 1936 analysis of computable functions\footnote{\cite{Tu36}.}.

Specifically, I consider an interpretation \(\mathcal{I}_{PA(\beta rouwer/Turing)}\) of PA in which we interpret the non-logical constants as in the standard interpretation \(\mathcal{I}_{PA(Stan-}\) \(_{dard/Tarski)}\) of PA in the `usual' manner, but interpret the logical constants algorithmically in an `unusual', yet consistent, manner.

This is the interpretation \(\mathcal{I}_{PA(\beta rouwer/Turing)}\) of PA obtained if, in Tarski's inductive definitions (1)-(8) below\footnote{cf.\ \cite{Me64}, pp.50-52.}---of the satisfaction and truth of the formulas of PA under the standard interpretation \(\mathcal{I}_{PA(Standard/Tarski)}\) of PA---we apply Occam's razor and weaken (1) by replacing it with the algorithmic definition of Turing-satisfiability (1b) (where (1a) is to be viewed simply as an intermediate definition of a catalytic nature).

\subsection{Defining satisfiability algorithmically}

Thus, we consider the definitions:

\vspace{+1ex}
\noindent (1) \textbf{Tarskian satisfiability}: If \([A(t_{1}, \ldots, t_{n})]\) is an atomic well-formed formula of PA, \(A^{*}(t^{*}_{1}, \ldots, t^{*}_{n})\) the corresponding interpretation under \(\mathcal{I}_{PA(Standard/Tarski)}\), and the sequence of PA-terms \([a_{1}], \ldots, [a_{n}]\) interprets under \(\mathcal{I}_{PA(Standard/Tarski)}\) as \(a^{*}_{1} \ldots, a^{*}_{n}\), then \(a^{*}_{1} \ldots, a^{*}_{n}\) satisfies \([A(t_{1}, \ldots, t_{n})]\) under \(\mathcal{I}_{PA(Standard/Tarski)}\) if, and only if, \(A^{*}(a^{*}_{1}, \ldots, a^{*}_{n})\) holds in the domain of the natural numbers;

\vspace{+1ex}
\noindent (1a) \textbf{Markovian satisfiability}: If \([A(t_{1}, \ldots, t_{n})]\) is an atomic well-formed formula of PA, \(A^{*}(t^{*}_{1}, \ldots, t^{*}_{n})\) the corresponding interpretation under \(\mathcal{I}_{PA(Standard-}\) \(_{/Tarski)}\), and the sequence of PA-terms \([a_{1}], \ldots, [a_{n}]\) interprets under \(\mathcal{I}_{PA(Standa-}\) \(_{rd/Tarski)}\) as \(a^{*}_{1} \ldots, a^{*}_{n}\), then \(a^{*}_{1} \ldots, a^{*}_{n}\) satisfies \([A(t_{1}, \ldots, t_{n})]\) under \(\mathcal{I}_{PA(Markov)}\) if, and only if, there is a Markov algorithm\footnote{See \cite{Me64}, p.209; \cite{Ma54}.} to establish that \(A^{*}(a^{*}_{1}, \ldots, a^{*}_{n})\) holds in the domain of the natural numbers;
\vspace{+1ex}

\noindent (1b) \textbf{Brouwer-Turing satisfiability}: If \([A(t_{1}, \ldots, t_{n})]\) is an atomic well-formed formula of PA, \(A^{*}(t^{*}_{1}, \ldots, t^{*}_{n})\) the corresponding interpretation under \(\mathcal{I}_{PA(Standa-}\) \(_{rd/Tarski)}\), and the sequence of PA-terms \([a_{1}], \ldots, [a_{n}]\) interprets under \(\mathcal{I}_{PA(Stan-}\) \(_{dard/Tarski)}\) as \(a^{*}_{1} \ldots, a^{*}_{n}\), then \(a^{*}_{1} \ldots, a^{*}_{n}\) satisfies \([A(t_{1}, \ldots, t_{n})]\) under \(\mathcal{I}_{PA(\beta rouwer}\) \(_{/Turing)}\) if, and only if, \(A^{*}(a^{*}_{1}, \ldots, a^{*}_{n})\) is TM\(_{A^{*}}\) -computable as true (when treated as a Boolean expression) in the domain of the natural numbers, where TM\(_{A^{*}}\) is the Turing-machine defined by \(A^{*}(t^{*}_{1}, \ldots, t^{*}_{n})\).

\subsection{Defining truth under an interpretation}

Further, both \(\mathcal{I}_{PA(Standard/Tarski)}\) and \(\mathcal{I}_{PA(\beta rouwer/Turing)}\) identically define the satisfaction and truth of the compound formulas of PA inductively as usual under the corresponding interpretation  as follows:

\vspace{+1ex}
\noindent (2) A sequence \(s\) satisfies \([\neg A]\) under \(\mathcal{I}_{PA(Standard/Tarski)}/\mathcal{I}_{PA(\beta rouwer/Turing)}\) if, and only if, \(s\) does not satisfy \([A]\);

\vspace{+1ex}
\noindent (3) A sequence \(s\) satisfies \([A \rightarrow B]\) under \(\mathcal{I}_{PA(Standard/Tarski)}/\mathcal{I}_{PA(\beta rouwer/Turing)}\) if, and only if, either \(s\) does not satisfy \([A]\), or \(s\) satisfies \(B\);

\vspace{+1ex}
\noindent (4) A sequence \(s\) satisfies \([(\forall x_{i})A]\) under \(\mathcal{I}_{PA(Standard/Tarski)}/\mathcal{I}_{PA(\beta rouwer/Turing)}\) if, and only if, every denumerable sequence in the domain \(D\) of \(\mathcal{I}_{PA(Standard/Tarski)}\) \(/\mathcal{I}_{PA(\beta rouwer/Turing)}\) which differs from \(s\) in at most the \(i\)'th component satisfies \(A\);

\vspace{+1ex}
\noindent (5) A well-formed formula \(A\) of PA is true under \(\mathcal{I}_{PA(Standard/Tarski)}/\mathcal{I}_{PA(\beta rou-}\) \(_{wer}\) \(_{/Turing)}\) if, and only if, every denumerable sequence in the domain \(D\) of \(\mathcal{I}_{PA(Standard/Tarski)}/\mathcal{I}_{PA(\beta rouwer/Turing)}\) satisfies \(A\);

\vspace{+1ex}
\noindent (6) A well-formed formula \(A\) of PA is false under \(\mathcal{I}_{PA(Standard/Tarski)}/\mathcal{I}_{PA(\beta ro-}\) \(_{uwer/Turing)}\) if, and only if, no sequence in the domain \(D\) of \(\mathcal{I}_{PA(Standard/Tarski)}/\) \(\mathcal{I}_{PA(\beta rouwer/Turing)}\) satisfies \(A\).

\subsection{Interpreting the universal quantifier effectively}

It follows from these definitions that:

\vspace{+1ex}
\noindent (7) \textbf{Tarskian universality}: A PA-formula such as \([(\forall x)A(x)]\) interprets as true under \(\mathcal{I}_{PA(Standard/Tarski)}\) if, and only if, for any given natural number \(n\), \(A^{*}(n)\) is true, where \(A^{*}(x)\) is the interpretation of \([A(x)]\) under \(\mathcal{I}_{PA(Standard/Tarski)}\);

\vspace{+1ex}
\noindent (7a) \textbf{Markovian universality}: A PA-formula such as \([(\forall x)A(x)]\) interprets as true under \(\mathcal{I}_{PA(Markov)}\) if, and only if, there is a \textit{Markov algorithm} that, for any given natural number \(n\), will establish \(A^{*}(n)\) as true, where \(A^{*}(x)\) is the interpretation of \([A(x)]\) under \(\mathcal{I}_{PA(Standard/Tarski)}\);

\vspace{+1ex}
\noindent (7b) \textbf{Brouwer-Turing universality}: A PA-formula such as \([(\forall x)A(x)]\) interprets as true under \(\mathcal{I}_{PA(\beta rouwer/Turing)}\) if, and only if, for any given natural number \(n\), the Turing-machine TM\(_{A{*}}\) computes \(A^{*}(n)\)---when treated as a Boolean function---as true, where \(A^{*}(x)\) is the interpretation of \([A(x)]\) under \(\mathcal{I}_{PA(Standard/Tarski)}\), and TM\(_{A{*}}\) is the Turing-machine defined by the quantifier-free expression in the prenex normal form of \(A^{*}\).

\subsection{Interpreting the existential quantifier effectively}

Further:

\vspace{+1ex}
\noindent (8) \textbf{Tarskian particularisation}: A PA-formula such as \([(\exists x)A(x)]\)---the abbreviation for \([\neg(\forall x)\neg A(x)]\)---interprets as true under \(\mathcal{I}_{PA(Standard/Tarski)}\) if, and only if, it is not true that, for any given natural number \(n\), \(A^{*}(n)\) is false, where \(A^{*}(x)\) is the interpretation of \([A(x)]\) under \(\mathcal{I}_{PA(Standard/Tarski)}\), \textit{and} we may conclude that there exists some natural number \(n\) such that \(A^{*}(n)\) holds;

\vspace{+1ex}
\noindent (8a) \textbf{Markovian particularisation}: A PA-formula such as \([(\exists x)A(x)]\)---the abbreviation for \([\neg(\forall x)\neg A(x)]\)---interprets as true under \(\mathcal{I}_{PA(Markov)}\) if, and only if, it is not true that there is an \textit{Markov algorithm} that, for any given natural number \(n\), will establish \(A^{*}(n)\) as false, where \(A^{*}(x)\) is the interpretation of \([A(x)]\) under \(\mathcal{I}_{PA(Standard/Tarski)}\), \textit{but} we may \textit{not} conclude that there exists some natural number \(n\) such that \(A^{*}(n)\) holds (since \(A^{*}(x)\) may be an instantiationally, but not algorithmically, decidable relation such that, for any given natural number \(n\), \(A^{*}(n)\) is false);

\vspace{+1ex}
\noindent (8b) \textbf{Brouwer-Turing particularisation}: A PA-formula such as \([(\exists x)A(x)]\)---the abbreviation of \([\neg(\forall x)\neg A(x)]\)---interprets as true under \(\mathcal{I}_{PA(\beta rouwer/Turing)}\) if, and only if, it is not true that, for any given natural number \(n\), the Turing-machine TM\(_{A^{*}}\) computes \(A^{*}(n)\)---when treated as a Boolean function---as false, where \(A^{*}(x)\) is the interpretation of \([A(x)]\) under \(\mathcal{I}_{PA(Standard/Tarski)}\), and TM\(_{A^{*}}\) is the Turing-machine defined by the quantifier-free expression in the prenex normal form of \(A^{*}\), \textit{but} we may \textit{not} conclude that there exists some natural number \(n\) such that \(A^{*}(n)\) holds (since \(A^{*}(x)\) may be a Halting-type of relation such that, for any given natural number \(n\), \(A^{*}(n)\) is false).

\subsection{The finitary interpretation \(\mathcal{I}_{PA(\beta rouwer/Turing)}\) of PA is sound}

Now, the PA-axioms PA\(_{1}\) to PA\(_{8}\)---which do not involve any quantification---interpret straightforwardly as Turing-computably always-true arithmetical relations under \(\mathcal{I}_{PA(\beta rouwer/Turing)}\).

Further, for any given PA-formula \([F(x)]\), the Induction axiom schema PA\(_{9}\) interprets under \(\mathcal{I}_{PA(\beta rouwer/Turing)}\) as a tautology, since:

\begin{quote}
`\([F(0) \vee (\forall x)(F(x) \rightarrow F(x'))]\) is true under \(\mathcal{I}_{PA(\beta rouwer/Turing)}\)'

`\([(\forall x)F(x)]\) is true under \(\mathcal{I}_{PA(\beta rouwer/Turing)}\)'
\end{quote}

\noindent both mean:

\begin{quote}
`For any given natural number \(n\), the Turing-machine TM\(_{F^{*}}\) computes \(F^{*}(n)\) as true'
\end{quote} 

\noindent where \(F^{*}(x)\) is the interpretation of \([F(x)]\) under \(\mathcal{I}_{PA(\beta rouwer/Turing)}\).

Similarly, Generalisation too interprets under \(\mathcal{I}_{PA(\beta rouwer/Turing)}\) as a tautology since, again:

\begin{quote}
`\([F(x)]\) is true (i.e., satisfied for any given \(x\)) under \(\mathcal{I}_{PA(\beta rouwer/Turing)}\)'

`\([(\forall x)F(x)]\) is true under \(\mathcal{I}_{PA(\beta rouwer/Turing)}\)'
\end{quote}

\noindent both mean:

\begin{quote}
`For any given natural number \(n\), the Turing-machine TM\(_{F^{*}}\) computes \(F^{*}(n)\) as true'
\end{quote} 

\noindent where \(F^{*}(x)\) is the interpretation of \([F(x)]\) under \(\mathcal{I}_{PA(\beta rouwer/Turing)}\).

It is also straightforward to show that Modus Ponens preserves truth under \(\mathcal{I}_{PA(\beta rouwer/Turing)}\).

Thus the axioms of PA are constructively satisfied/true under the finitary interpretation \(\mathcal{I}_{PA(\beta rouwer/Turing)}\), and the rules of inference of PA preserve the properties of satisfaction/truth under \(\mathcal{I}_{PA(\beta rouwer/Turing)}\).

It follows that:

\begin{theorem}
: If a formula \([F]\) is a theorem of a first-order Peano Arithmetic then there is a Turing-machine TM\(_{F}\) that will compute the arithmetical proposition---or relation---\(F\) as true---or \textit{always} true (i.e., true for any natural number values assigned to the free variables of \(F\)), respectively---in a finite number of steps.
\end{theorem}

Hence the finitary interpretation \(\mathcal{I}_{PA(\beta rouwer/Turing)}\) of PA is sound, and so it defines a finitary model of PA.

\subsubsection{\(\mathcal{I}_{PA(\beta rouwer/Turing)}\) settles the Poincar\'{e}-Hilbert debate}

Moreover, since \(\mathcal{I}_{PA(\beta rouwer/Turing)}\) highlights the tautological character of both the Generalisation rule of inference and the Induction Schema in PA, it settles the Poincar\'{e}-Hilbert debate in the latter's favour. Poincar\'{e} maintained that the consistency of the method of induction could never be proven except through the inductive method itself; Hilbert believed a finitary proof of the consistency of PA was possible\footnote{See \cite{Hi27}, p.472; also \cite{We27}, p482; \cite{Br13}, p.59.}.

The above analysis suggests that their difference resolves in Hilbert's favour once we explicitly differentiate between `For all \ldots'---which implicitly, but unintendedly, implies algorithmic (and, ipso facto, instantiational) verifiability---and `For any given \ldots'---which implies only instantiational verifiability.

\subsection{The introduction of an algorithmic method for the decidability of satisfaction and truth under Tarski's definitions is necessary}

The question arises: Is the introduction of an algorithmic method for the decidability of satisfaction and truth under Tarski's definitions necessary?

I give an affirmative answer by defining a formal interpretation \(\mathcal{I}_{PA(G\ddot{o}del)}\) of PA, and show that---as Brouwer maintained---unrestricted interpretation of the universal quantifier can also lead to an inconsistency\footnote{\cite{Br08}; see also \cite{Br23}, p.336; \cite{Br27}, p.491; \cite{We27}, p.483.}.

I define satisfaction under \(\mathcal{I}_{PA(G\ddot{o}del)}\) as:

\vspace{+1ex}
\noindent (1c) \textbf{G\"{o}delian satisfiability}: If \([A(t_{1}, \ldots, t_{n})]\) is an atomic well-formed formula of PA then the sequence of PA-terms \([a_{1}], \ldots, [a_{n}]\) satisfies \([A(t_{1}, \ldots, t_{n})]\) under \(\mathcal{I}_{PA(G\ddot{o}del)}\) if, and only if, \([A(a_{1}, \ldots, a_{n})]\) is PA-provable.

\vspace{+1ex}
The satisfaction and truth of the compound formulas of PA are defined inductively under the interpretation \(\mathcal{I}_{PA(G\ddot{o}del)}\) as usual by (2) to (6) above.

If we accept that (1c) is a consistent definition of satisfiability, then it is straightforward to show that \(\mathcal{I}_{PA(G\ddot{o}del)}\) is sound. 

Now, it also follows from the definitions (1c) and (2) to (6) that:

\vspace{+1ex}
\noindent (7c) \textbf{G\"{o}delian universality}: A PA-formula such as \([(\forall x)A(x)]\) interprets as true under \(\mathcal{I}_{PA(G\ddot{o}del)}\) if, and only if, for any given numeral \([n]\), \([A(n)]\) is PA-provable;

\vspace{+1ex}
\noindent (8c) \textbf{G\"{o}delian particularisation}: A PA-formula such as \([(\exists x)A(x)]\)---the abbreviation for \([\neg(\forall x)\neg A(x)]\)---interprets as true under \(\mathcal{I}_{PA(G\ddot{o}del)}\) if, and only if, it is not true that for any given numeral \([n]\), \([\neg A(n)]\) is PA-provable.

\vspace{+1ex}
Further, we have shown above that we can construct a G\"{o}delian PA-formula \([R(x)]\) such that:

\vspace{+1ex}
(i) \([(\forall x)R(x)]\) is not PA-provable;

\vspace{+1ex}
(ii) \([\neg (\forall x)R(x)]\) is PA-provable;

\vspace{+1ex}
(iii) for any given numeral \([n]\), \([R(n)]\) is PA-provable.

\vspace{+1ex}
However, if \(\mathcal{I}_{PA(G\ddot{o}del)}\) is sound, then (ii) implies that it is not the case that, for any given numeral \([n]\), \([R(n)]\) is PA-provable. It follows that \(\mathcal{I}_{PA(G\ddot{o}del)}\), too, is not sound even though---unlike \(\mathcal{I}_{PA(Standard/Tarski)}\)---\(\mathcal{I}_{PA(G\ddot{o}del)}\) defines particularisation logically, and not platonically as Aristotle did (possibly influenced by his teacher Plato).

We conclude that if an interpretation of PA lacks specification of an effective method for determining the satisfaction of the atomic formulas of PA in definitions (1) and (1c), then this lack is reflected in the non-constructivity of definitions (4) to (6) under \(\mathcal{I}_{PA(Standard/Tarski)}\) and \(\mathcal{I}_{PA(G\ddot{o}del)}\); these are, then, not sufficient to distinguish between arithmetical relations that are algorithmically decidable, and those---such as G\"{o}del's \(R(x)\)---that are instantiationally, but not algorithmically, decidable.

\subsection{PA is categorical: A Provability Theorem for PA}

I now show that PA is categorical, and can have no non-standard model, since it is `algorithmically' complete in the sense that:

\begin{theorem}
(Provability Theorem for PA): A total arithmetical relation \(F(x)\)---when treated as a Boolean function---defines a Turing-machine TM\(_{F}\) which computes \(F(x)\) as \textit{always} true (i.e., true for any given natural number input \(n\)) if, and only if, the corresponding PA-formula [\(F(x)\)] is PA-provable.
\end{theorem}

\vspace{+1ex}
\noindent \textbf{Proof}: It follows from Turing's seminal 1936 paper \footnote{\cite{Tu36}.} that every quantifier-free arithmetical relation \(F(x)\) (when interpreted as a Boolean function) defines a Turing-machine TM\(_{F}\)\footnote{In the general case, TM\(_{F}\) is defined by the quantifier-free expression in the prenex normal form of \(F(x)\); cf.\ \cite{Tu36}, pp.\ 138-139} such that \(F(x)\) is TM\(_{F}\)-computable if, and only if, for \textit{any} given natural number \(n\), TM\(_{F}\) will compute \(F(n)\) as either true, or as false, over the structure \(\mathcal{N}\).

Now, we have that:

\begin{tabbing}
\indent (a) \= [\((\forall x)F(x)\)] is defined as true under the interpretation \(\mathcal{I}_{PA(\beta rouwer/Turing)}\) \\ \> if, and only if, the Turing-machine TM\(_{F}\) computes \(F(n)\) as \textit{always} true \\ \> (i.e., as true for any given natural number \(n\)) under \(\mathcal{I}_{PA(\beta rouwer/Turing)}\); \\

\indent (b) \> [\((\exists x)F(x)\)] is an abbreviation of [\(\neg (\forall x)\neg F(x)\)], and is defined as true \\ \> under \(\mathcal{I}_{PA(\beta rouwer/Turing)}\) if, and only if, it is not the case that the Turi- \\ \> ng-machine TM\(_{F}\) computes \(F(n)\) as \textit{always} false (i.e., as false for any \\ \> given natural number \(n\)) under \(\mathcal{I}_{PA(\beta rouwer/Turing)}\).
\end{tabbing}

Moreover, since \(\mathcal{I}_{PA(\beta rouwer/Turing)}\) is sound, it defines a finitary model of PA over \(\mathcal{N}\)---say \(\mathcal{M}_{PA(\beta)}\)---such that:

\vspace{+1ex}
If [\((\forall x)F(x)\)] is PA-provable, then the Turing-machine TM\(_{F}\) computes \(F(n)\) as \textit{always} true (i.e., as true for any given natural number \(n\)) in \(\mathcal{M}_{PA(\beta)}\);

If [\(\neg(\forall x)F(x)\)] is PA-provable, then it is not the case that the Turing-machine TM\(_{F}\) computes \(F(n)\) as \textit{always} true (i.e., as true for any given natural number \(n\)) in \(\mathcal{M}_{PA(\beta)}\).

\vspace{+1ex}
Further, we cannot have that both [\((\forall x)F(x)\)] and [\(\neg(\forall x)F(x)\)] are PA-unprovable for some PA-formula \(F(x)\), as this would yield the contradiction:

\begin{tabbing}
\indent (i) \= There is a finitary model---say \(M1_{\beta}\)---of PA+[\((\forall x)F(x)\)] 
in which the \\ \> Turing-machine TM\(_{F}\) computes \(F(n)\) as \textit{always} true (i.e., as true for \\ \> any given natural number \(n\)). \\

\indent (ii) \= There is a finitary model---say \(M2_{\beta}\)---of PA+[\(\neg(\forall x)F(x)\)] in which it is \\ \> not the case that the Turing-machine TM\(_{F}\) computes \(F(n)\) as \textit{always} \\ \> true (i.e., as true for any given natural number \(n\)).
\end{tabbing}

Hence PA is `algorithmically' complete, in the sense that a total arithmetical relation \(F(x)\)---when treated as a Boolean function---defines a Turing-machine TM\(_{F}\) which computes \(F(x)\) as \textit{always} true if, and only if, the corresponding PA-formula [\(F(x)\)] is PA-provable\footnote{Note that [\((\forall x)F(x)\)] is PA-provable if, and only if, [\(F(x)\)] is PA-provable.} \(\Box\)

Hence PA can have no `non-standard' model, and we have that:

\begin{corollary}
: PA is categorical.
\end{corollary}

\subsubsection{The significance of the Provability Theorem for PA}

\vspace{+1ex}
The Provability Theorem for PA is of particular interest computationally, since it formally expresses some implicitly held beliefs in interpretations of computability theory. For instance, in their paper `Passages of Proof'\footnote{\cite{CCS01}, p.13.}, the authors hold that:

\begin{quote}
Classically, there are two equivalent ways to look at the mathematical notion of proof: logical, as a finite sequence of sentences strictly obeying some axioms and inference rules, and computational, as a specific type of computation. Indeed, from a proof given as a sequence of sentences one can easily construct a Turing-machine producing that sequence as the result of some finite computation and, conversely, given a machine computing a proof we can just print all sentences produced during the computation and arrange them into a sequence.
\end{quote}

In other words, the authors seem to hold---echoing the sense of the Provability Theorem for PA---that Turing-computability of a `proof', in the case of an arithmetical proposition, is equivalent to provability of its representation in PA.

\subsubsection{G\"{o}del's arithmetical relation \(R(n)\) is effectively computable as true for any given natural number \(n\), but it is not Turing-computable as true for any given natural number \(n\)}

Another significant consequence is that the Provability Theorem for PA provides an intuitively plausible explanation for the fact that---following G\"{o}del---we can define an arithmetical relation \(R(n)\) that is true for any given natural number \(n\), but the corresponding PA-formula [\(R(x)\)] is not PA-provable\footnote{\cite{Go31}, p25(1).}.

\vspace{+1ex}
(By Generalisation\footnote{\textit{Generalisation in PA}: [\((\forall x)F\)] follows from [\(F\)].}, stating that the PA-formula [\(R(x)\)] is not PA-provable is equivalent to stating that the PA-formula [\((\forall x)R(x)\)]\footnote{G\"{o}del defines, and refers to, the formula corresponding to this in his formal system P by its G\"{o}del-number \(17Genr\).} is not PA-provable; the latter is what G\"{o}del actually proved for his formal system of Peano Arithmetic P.)

\vspace{+1ex}
It simply means that, given any natural number \(n\), \(R(n)\) is true but---as a consequence of the Provability Theorem for PA---the Turing-machine TM\(_{R}\) does \textit{not} compute \(R(n)\) as \textit{always} true (i.e., true for any given natural number \(n\)).

The counter-part of this statement in PA can be expressed as:

\begin{quote}
Given any PA-numeral \([n]\), there is always some PA-deduction \(D_{n}\)---i.e., a finite proof-sequence of PA-formulas---whose last formula is \([R(n)]\), but there is no common deduction \(D\) which is independent of \([n]\), i.e., there is no proof-sequence of PA-formulas whose last formula is \([R(x)]\), so \([R(x)]\) is PA-unprovable.
\end{quote}

In other words we have as a consequence of the Provability Theorem for PA that:

\begin{theorem}
: (First Tautology Theorem) The Turing-machine TM\(_{R}\) does not compute the tautology \(R(n)\)---when treated as a Boolean function---as \textit{always} true (i.e., true for any given natural number \(n\)).
\end{theorem}

\textit{Proof}: In his seminal 1931 paper\footnote{\cite{Go31}.}, G\"{o}del has constructed an arithmetical relation \(R(n)\) that is meta-mathematically provable as true for any given natural number \(n\) but --- since the corresponding PA-formula [\(R(x)\)]\footnote{G\"{o}del defines, and refers to, this formula by its G\"{o}del-number \(r\) (cf.\ [Go31], p25, eqn.12).} is not PA-provable\footnote{\cite{Go31}, p25(1).}---it follows from the Provability Theorem for PA that there is no Turing-machine TM\(_{R}\) that computes \(R(n)\) as \textit{always} true (i.e. true for any given natural number \(n\)).\(\Box\)

\subsubsection{P\(\neq\)NP}

The First Tautology Theorem has a curious consequence concerning the celebrated PvNP problem\footnote{\cite{Cook}.}.

In a paper presented to ICM 2002, Ran Raz comments\footnote{\cite{Ra02}.}:

\begin{quote}
``A Boolean formula \(f(x_{1}, \ldots, x_{n})\) is a tautology if \(f(x_{1}, \ldots, x_{n}) = 1\) for every \(x_{1}, \ldots, x_{n}\). A Boolean formula \(f(x_{1}, \ldots, x_{n})\) is unsatisfiable if \(f(x_{1}, \ldots, x_{n}) = 0\) for every \(x_{1}, \dots, x_{n}\). Obviously, \(f\) is a tautology if and only if \(\neg f\) is unsatisfiable.

Given a formula \(f(x_{1}, \ldots, x_{n})\), one can decide whether or not \(f\) is a tautology by checking all the possibilities for assignments to \(x_{1}, \ldots, x_{n}\). However, the time needed for this procedure is exponential in the number of variables, and hence may be exponential in the length of the formula \(f\). \ldots

P\(\neq\)NP is the central open problem in complexity theory and one of the most important open problems in mathematics today. The problem has thousands of equivalent formulations. One of these formulations is the following:

Is there a polynomial time algorithm \(\mathcal{A}\) that gets as input a Boolean formula \(f\) and outputs 1 if and only if \(f\) is a tautology?

P\(\neq\)NP states that there is no such algorithm."
\end{quote}

So, if \(R(n)\) is constructively computable as a tautology, but not recognisable as a tautology by any Turing-machine, then, by the First Tautology Theorem:

\begin{corollary}
: P\(\neq\)NP!
\end{corollary}

\subsubsection{Is the solution P\(\neq\)NP significant?}

However, the PvNP problem assumes the significance usually accorded to it only to the extent that its solution throws light on the practical consequences---mentioned above by Raz---that follow from the computational complexity of a number-theoretic relation (or function), and not simply from the philosophical consequences of its logical properties.

The following shows why a trivial logical solution of the PvNP problem---such as that indicated above, which addresses the problem as it is presently formulated---may not be significant from a computational complexity perspective.

\subsubsection{Can the PvNP problem be formulated to highlight the significance of computational complexity?}

In his 1931 paper\footnote{\cite{Go31}.}, G\"{o}del specifically defined his arithmetical relation \(R(n)\) so that it is instantiationally equivalent to a primitive recursive relation, \(Q(n)\), which, of course, is Turing-computable as true for any given natural number \(n\).

So, even if P\(\neq\)NP because PA has no non-standard models, we still have two number-theoretic relations that are instantiationally equivalent even though they are not algorithmically (mathematically) identical.

The equivalence should, prima facie, suffice to define the computational complexity involved in one case as notionally equal to that actually involved in the other.

\subsubsection{\textit{Why} the PvNP problem may \textit{remain} intractable in ZF}

Now, the above distinction between instantiational equivalence and mathematical identity would be absent in any set-theoretic approach to the PvNP problem (and hence in its usual set-theoretic interpretation in terms of the classes of decidable and recursively enumerable languages), since it is an axiom of ZF that two ZF relations (or functions) are instantiationally equivalent if, and only if, they are set-theoretically (\textit{mathematically}?) identical.

In other words, if we accept that the consistency of PA implies that P\(\neq\)NP, then, if ZF is consistent, it can neither reflect that P\(\neq\)NP nor that P\(=\)NP.

This suggests that it may be better to reformulate the PvNP problem number-theoretically so as to avoid trivial \textit{logical} solutions in PA, or \textit{logical} blind alleys in ZF, and to reflect better its computational significance.

\subsubsection{The standard Church and Turing Theses}

Another consequence of the Provability Theorem for PA is that:

\begin{theorem}
: (Second Tautology Theorem) G\"{o}del's tautology \(R(n)\) is constructively computable as \textit{always} true (i.e., true for any given natural number \(n\)).
\end{theorem}

\textit{Proof}: G\"{o}del has defined a primitive recursive relation, \(xB_{PA}y\) that holds if, and only if, \(y\) is the G\"{o}del-number of a PA-formula, say [\(Y\)], and \(x\) the G\"{o}del-number of a PA-proof of [\(Y\)] ([Go31], p22, dfn.\ 45).

Since every primitive recursive relation is Turing-computable (when treated as a Boolean function), \(xB_{PA}y\) defines a Turing-machine TM\(_{B}\) that halts on any natural number values of \(x\) and \(y\).

Now, if \(g_{[R(1)]}\), \(g_{[R(2)]}\), \ldots are the G\"{o}del-numbers of the PA-formulas [\(R(1)\)], [\(R(2)\)], \dots, it follows that, for any given natural number \(n\), when  the natural number value \(g_{[R(n)]}\) is input for \(y\), the Turing-machine TM\(_{B}\) must halt for some value of \(x\)---which is the G\"{o}del-number of some PA-proof of [\(R(n)\)]---since G\"{o}del has shown\footnote{\cite{Go31}, p25(1).} that [\(R(n)\)] is PA-provable for any given numeral [\(n\)].

Hence \(R(n)\) is constructively computable as true for any given natural number \(n\). \(\Box\)

The Second Tautology Theorem also has an interesting Corollary.

Now, we can show that:

\begin{quote}
(a)	Every Turing-computable function (or relation, treated as a Boolean function) \(F\) is partial recursive , and, if \(F\) is total , then \(F\) is recursive\footnote{cf.\ \cite{Me64}, p.233, Corollary 5.13.}.

(b)	Every partial recursive function (or relation, treated as a Boolean function) is Turing-computable\footnote{cf.\ \cite{Me64}, p.237, Corollary 5.15.}.
\end{quote}

It follows that the following---essentially unverifiable but refutable---theses are equivalent\footnote{cf.\ \cite{Me64}, p.237.}:

\begin{quote}
\textbf{Standard Church's Thesis}\footnote{\textit{Church's (original) Thesis}: The effectively computable number-theoretic functions are the algorithmically computable number-theoretic functions \cite{Ch36}.}: A number-theoretic function (or relation, treated as a Boolean function) is effectively computable if, and only if, it is partial-recursive\footnote{cf.\ \cite{Me64}, p.227.}.

\textbf{Standard Turing's Thesis}\footnote{After describing what he meant by ``computable" numbers in the opening sentence of his seminal paper on Computable Numbers \cite{Tu36}, Turing immediately expressed this thesis---albeit informally---as: ``\ldots the computable numbers include all numbers which could naturally be regarded as computable".}: A number-theoretic function (or relation, treated as a Boolean function) is effectively computable if, and only if, it is Turing-computable\footnote{cf.\ \cite{BBJ03}, p.33.}.
\end{quote}

However, it follows from the Second Tautology Theorem that, since there is a number-theoretic relation which---treated as a Boolean function---is (effectively) constructively computable but not (algorithmically) Turing-computable:

\begin{theorem}
: The Church and Turing theses are false.
\end{theorem}

\subsubsection{Church's Thesis as a weaker equivalence}

However, instead of expressing it strongly as a refutable identity, Church's Thesis can also be postulated as the weaker---and intuitively more plausible---equivalence:

\begin{quote}
\textbf{Weak Church's Thesis}: A number-theoretic function (or relation, treated as a Boolean function) is effectively computable if, and only if, it is instantiationally equivalent to a partial-recursive function.
\end{quote}

It is significant that G\"{o}del (initially) and Church (subsequently - not least because of G\"{o}del's disquietitude) enunciated Church's formulation of `effective computability' as a Thesis because G\"{o}del was instinctively uncomfortable with accepting it as a definition that fully captures the essence of `\textit{intuitive} effective computability'\footnote{cf.\ \cite{Si97}}.

G\"{o}del's reservations seem vindicated if we accept that a number-theoretic function can be effectively computable instantiationally, but not algorithmically.

Since every algorithmically computable function is, necessarily, computable instantiationally, we can now define:

\begin{quote}
\textbf{Definition}: A number-theoretic function is effectively computable (\textit{intuitively}) if, and only if, it is effectively computable instantiationally.
\end{quote}

We thus see that, in its standard formulations, Church's Thesis violates the doctrine of Occam's razor by postulating a strong identity---and not simply a weak equivalence---between an effectively computable number-theoretic function and some partial recursive function.

Consequently, the Thesis does not admit the possibility of number-theoretic relations that are constructively computable (i.e., effectively computable instantiationally), but not Turing-computable (i.e., not effectively computable algorithmically).

\section{Cantor and ZF set theory}

In this section I consider the belief that all \textit{significant} mathematical `truths' --- such as, for example, the theorems of a first-order Peano Arithmetic (PA) --- can be interpreted as theorems of a set-theory such as ZF without \textit{any} loss of generality.

For instance, in the chapter, ``NP and NP completeness", of their forthcoming book, authors Arora and Barak write that\footnote{\cite{Ar08}, p2.24(60), Ex.6, Ch.2.}:

\begin{quote}
Mathematics can be axiomatized using for example the Zermelo Frankel system, which has a finite description.
\end{quote}

That such a belief is almost universal today is a reflection of the fact that---for over a generation---it has been explicitly echoed in standard texts with increasing certitude:

\begin{quote}
\ldots NBG\footnote{Also `ZF'.} apparently can serve as a foundation for all present-day mathematics (i.e., it is clear to every mathematician that every mathematical theorem can be translated and proved within NBG, or within extensions of NBG obtained by adding various extra axioms such as the Axiom of Choice) \ldots \\ \ldots \textit{Mendelson\footnote{\cite{Me64}, p193.}}

Such is the case, for instance, with the formal systems considered in works on set theory, such as the one known as ZFC, which are adequate for formalizing essentially all accepted mathematical proofs. \\ \ldots \textit{Boolos, Burgess, and Jeffrey\footnote{\cite{BBJ03}, p225.}}

\end{quote}

\subsection{Relativising PA in ZF}

Now a standard method of interpreting PA in ZF is by relativising PA in ZF\footnote{By appealing to Cantor's first limit ordinal, \(\omega\).} as follows\footnote{cf.\ \cite{Me64}, p192.}:

\begin{quote}

Given any wf \(\mathcal{A}\) of formal number theory S\footnote{Read as `PA'.} \ldots, we can associate with \(\mathcal{A}\) a wf \(\mathcal{A}^*\) of NBG\footnote{Read as `ZF'.} as follows: first replace every ``+" by ``\(+_{0}\)", and every ``."\footnote{Read as ``\(\star\)".} by ``\(\times_{0}\)"\footnote{Note that ``\(+_{0}\)" and ``\(\times_{0}\)" denote ordinal addition and ordinal multiplication.}; then, if \(\mathcal{A}\) is \(\mathcal{B} \supset \mathcal{C}\), or \(\neg \mathcal{B}\), respectively, and we already have found \(\mathcal{B}^*\) and \(\mathcal{C}^*\), let \(\mathcal{A}^*\) be \(\mathcal{B}^* \supset \mathcal{C}^*\), or \(\neg (\mathcal{B}^*)\), respectively; if \(\mathcal{A}\) is \((x\footnote{Read as `\(\forall x\)'.})\mathcal{B} (x)\), replace it by \((x)(x\in \omega \supset \mathcal{B}^*(x)\). This completes the definition of \(\mathcal{A}^*\).

Now, if \(x_{1}, \ldots, x_{n}\) are the free variables of \(\mathcal{A}\), prefix `\((x_{1} \supset \omega \wedge \ldots \wedge x_{n} \supset \omega)\supset\)' to \(\mathcal{A}^*\), obtaining a wf \(\mathcal{A}\#\). This amounts to restricting all variables to \(\omega\) and interpreting addition, multiplication, and the successor function on integers as the corresponding operations on ordinals\footnote{Note that `\(x'\)' interprets as `\(x\cup {x}\)'.}.

Then every axiom \(\mathcal{A}\) of S is transformed into a theorem \(\mathcal{A}\#\) of NBG.

\ldots Applications of modus ponens are preserved under the transformations of \(\mathcal{A}\) into \(\mathcal{A}\#\). Also \ldots applications of the Generalisation Rule lead from theorems to theorems.

Therefore, for every theorem \(\mathcal{A}\) of S, \(\mathcal{A}\#\) is a theorem of NBG, and we can translate into NBG all the theorems of S \ldots

\end{quote}

\subsection{\textit{Can} ZF yield all significant mathematical `truths'?}

However, in a 1992 talk\footnote{\cite{Fe92}.}, Feferman tacitly sounded a cautionary note on an unqualified understanding --- such as that cited above --- of the ambit of set-theory.

\begin{quote}
\ldots... we are mainly concerned with the relation \(T_{1} \leq T_{2}\) when [\textit{theory}] \(T_{1}\) is a part of \(T_{2}\), either directly or by translation. In contrast, \(T_{2}\) tends to be more comprehensive than \(T_{1}\) in the case of interpretations; a familiar example is Peano Arithmetic (as \(T_{1}\)) in Zermelo-Fraenkel ZF (as \(T_{2}\)), where the natural numbers are interpreted as the finite ordinals. This is a conceptual reduction of number theory to set theory, but not a foundational reduction, because the latter system is justified only by an uncountable infinitary framework whereas the former is justified simply by a countable infinitary framework.
\end{quote}

A significant point which emerges from Feferman's talk is that we may not appeal unrestrictedly to set-theoretical reasoning when studying the foundational framework of PA. As we show below, the cautionary element underlying Feferman's remarks is particularly relevant when applied to two issues that are primarily rooted in number-theoretic --- and not set-theoretic --- concerns.

\subsection{Why PA \textit{cannot} admit a set-theoretical model}

Let [\(G(x)\)] denote the PA-formula:

\addvspace{+1ex}
\([x=0 \vee \neg(\forall y)\neg(x=y^{\prime})]\)

\addvspace{+1ex}
This translates, under every unrelativised interpretation of PA, as:

\addvspace{+1ex}
If \(x\) denotes an element in the domain of an unrelativised interpretation of PA, either \(x\) is 0, or \(x\) is a `successor'.

\addvspace{+1ex}
Further, in every such interpretation of PA, if \(G(x)\) denotes the interpretation of [\(G(x)\)]:

\addvspace{+1ex}
(a)	\(G(0)\) is true;

(b)	If \(G(x)\) is true, then \(G(x^{\prime})\) is true.

\addvspace{+1ex}
Hence, by G\"{o}del's completeness theorem:

\addvspace{+1ex}
(c)	PA proves \([G(0)]\);

(d)	PA proves \([G(x) \rightarrow G(x^{\prime})]\).

\addvspace{+1ex}
Further, by Generalisation:

\addvspace{+1ex}
(e)	PA proves \([(\forall x)(G(x) \rightarrow G(x^{\prime}))]\);

\addvspace{+1ex}
Hence, by Induction:

\addvspace{+1ex}
(f)	\([(\forall x)G(x)]\) is provable in PA.

\addvspace{+1ex}
In other words, except 0, every element in the domain of any unrelativised interpretation of PA is a `successor'. Further, \(x\) can only be a `successor' of a unique element in any such interpretation of PA. 

\subsubsection{PA and ZF have no common model}

Now, since Cantor's first limit ordinal, \(\omega\), is not the `successor' of any ordinal in the sense required by the PA axioms, and if there are no infinitely descending sequences of ordinals\footnote{\cite{Me64}, p261.} in a model---if any---of set-theory, PA and Ordinal Arithmetic\footnote{\cite{Me64}, p.187.} cannot have a common model, and so we cannot consistently extend PA to ZF simply by the addition of more axioms.

\subsubsection{Why the usual argument \textit{for} a non-standard model of PA is unconvincing}

Further, although we \textit{can} define a model of Arithmetic with an infinite descending sequence of elements\footnote{eg.\ \cite{BBJ03}, Section 25.\ 1, p303.}, any such model would be isomorphic to the ``\textit{true arithmetic}"\footnote{\cite{BBJ03}.\ p150.\ Ex.\ 12.\ 9.} of the integers (\textit{negative plus positive}), and \textit{not} to any model of PA\footnote{\cite{BBJ03}.\ Corollary 25.\ 3, p306.}.

Moreover---as we show in the next section---we cannot assume that we can consistently add a constant \(c\) to PA, along with the denumerable axioms [\(\neg (c = 0)\)], [\(\neg (c = 1)\)], [\(\neg (c = 2)\)], \ldots, since this would presume that which is sought to be proven, viz., that PA has a non-standard model.

We \textit{cannot} therefore---as suggested in standard texts\footnote{eg.\ \cite{BBJ03}.\ p306; \cite{Me64}, p112, Ex.\ 2.}---apply the Compactness Theorem and the (\textit{upward}) L\"{o}wenheim-Skolem Theorem to conclude that PA has a non-standard model!

\begin{quote}
\footnotesize{
\noindent \textbf{Compactness Theorem}: If every finite subset of a set of sentences has a model, then the whole set has a model\footnote{\cite{BBJ03}.\ p147.}.

\noindent \textbf{Upward L\"{o}wenheim-Skolem Theorem}: Any set of sentences that has an infinite model has a non-denumerable model\footnote{\cite{BBJ03}.\ p163.}.
}
\end{quote}

\subsection{A formal argument \textit{for} a non-standard model of PA}

The following argument\footnote{\cite{Lu08}.} attempts to validate the above line of reasoning suggested by standard texts \textit{for} the existence of non-standard models of PA:

\begin{tabbing}
\indent 1. \= Let \(<\)\(N\) (\textit{the set of natural numbers}); \(=\) (\textit{equality}); \('\) (\textit{the successor fun-} \\ \> \textit{ction}); \(+\) (\textit{the addition function}); \( \ast \) (\textit{the product function}); \(0\) (\textit{the null} \\ \> \textit{element})\(>\) be the structure, say [\(N\)], that serves to define a \textit{sound} int- \\ \> erpretation of PA. \\

\indent 2. \> Let T[\(N\)] be the set of PA-formulas that are satisfied or true in [\(N\)]. \\

\indent 3. \> The PA-provable formulas form a subset of T[\(N\)]. \\

\indent 4. \> Let \( \Gamma \) be the countable set of all PA-formulas of the form \([c_{n} = (c_{n+1})']\), \\ \> where the index \(n\) is a natural number. \\

\indent 5. \> Let T be the union of \( \Gamma \) and T[\(N\)]. \\

\indent 6. \> T[\(N\)] plus any finite set of members of \( \Gamma \) has a model, e.g., [\(N\)] itself, \\ \> since [\(N\)] is a model of any finite descending chain of successors. \\ 

\indent 7. \> Consequently, by Compactness, T has a model; call it \(M\). \\

\indent 8. \> \(M\) has an infinite descending sequence with respect to \('\) because it is a \\ \> model of \( \Gamma \). \\

\indent 9. \> Since PA is a subset of T, \(M\) is a non-standard model of PA.
\end{tabbing}

Now, if---as claimed above---[\(N\)] is a model of T[\(N\)] plus any finite set of members of \( \Gamma \), then all PA-formulas of the form \([c_{n} = (c_{n+1})']\) are PA-provable, \( \Gamma \) is a proper sub-set of the PA-provable formulas, and T is identically  T[\(N\)].

(The argument cannot be that some PA-formula of the form \([c_{n} = (c_{n+1})']\) is true in [\(N\)], but not PA-provable, as this would imply that PA+\([\neg (c_{n} = (c_{n+1})')]\) has a model other than [\(N\)]; in other words, it would presume that PA has a non-standard model.\footnote{The same objection applies to the usual argument found in standard texts (eg.\ [BBJ03].\ p306; \cite{Me64}, p112, Ex.\ 2.) which, again, is essentially that, if PA has a non-standard model at all, then one such model is obtained by assuming we can consistently add a single non-numeral constant \(c\) to the language of PA, and the countable axioms \(c \neq 0\), \(c \neq 1\), \(c \neq 2\), \ldots to PA. However, as noted earlier, this argument too does not resolve the question of whether such assumption validly allows us to conclude that there is a non-standard model of PA in the first place.

To place this distinction in perspective, Legendre and Gauss independently conjectured in 1796 that, if \(\pi (x)\) denotes the number of primes less than \(x\), then \(\pi (x)\) is asymptotically equivalent to \(x\)/In\((x)\). Between 1848/1850, Chebyshev confirmed that if \(\pi (x)\)/\{\(x\)/In\((x)\)\} has a limit, then it must be 1. However, the crucial question of whether \(\pi (x)\)/\{\(x\)/In\((x)\)\} has a limit at all was answered in the affirmative independently by Hadamard and de la Vall\'{e}e Poussin only in 1896.})

Consequently, the postulated model \(M\) of T in (7), by ``Compactness", is the model [\(N\)] that defines T[\(N\)]. However, [\(N\)] has no infinite descending sequence with respect to the successor function \('\), even though it is a model of \( \Gamma \). Hence the argument does not establish the existence of a non-standard model of PA with an infinite descending sequence with respect to the successor function \('\).

\subsection{The (\textit{upward}) Skolem-L\"{o}wenheim theorem applies only to first-order theories that admit an axiom of infinity}

We note, moreover, that the non-existence of non-standard models of PA would not contradict the (\textit{upward}) Skolem-L\"{o}wenheim theorem, since the proof of this theorem implicitly limits its applicability amongst first-order theories to those that are consistent with an axiom of infinity---in the sense that the proof implicitly requires that a constant, say \(c\), along with a denumerable set of axioms to the effect that \(c \neq 0, c \neq 1, \ldots \), can be consistently added to the theory. However, as seen in the previous section, this is not the case with PA.

\subsection{\textit{Why} PA has no set-theoretical model}

We can define the usual order relation `\(<\)' in PA so that every instance of the Induction Axiom schema, such as, say:

\vspace{+1ex}
(i) [\(F(0) \rightarrow ((\forall x)(F(x) \rightarrow F(x')) \rightarrow (\forall x)F(x))\)]

\vspace{+1ex}
yields the PA theorem:

\vspace{+1ex}
(ii) [\(F(0) \rightarrow ((\forall x) ((\forall y)(y < x \rightarrow F(y)) \rightarrow F(x)) \rightarrow (\forall x)F(x))\)]

\vspace{+1ex}
Now, if we interpret PA without relativisation in ZF in the sense indicated by Feferman\footnote{\cite{Fe92}.} --- i.e., numerals as finite ordinals, [\(x'\)] as [\(x \cup \left \{ x \right \}\)], etc. --- then (ii) always translates in ZF as a theorem:

\vspace{+1ex}
(iii) [\(F(0) \rightarrow ((\forall x)((\forall y)(y \in x \rightarrow F(y)) \rightarrow F(x)) \rightarrow (\forall x)F(x))\)]

\vspace{+1ex}
However, (i) does not always translate similarly as a ZF-theorem (which is why PA and ZF can have no common model), since the following is not necessarily provable in ZF:

\vspace{+1ex}
(iv) [\(F(0) \rightarrow ((\forall x)(F(x) \rightarrow F(x \cup \left \{x\right \})) \rightarrow (\forall x)F(x))\)]

\vspace{+1ex}
\textit{Example}: Define [\(F(x)\)] as `[\(x \in \omega\)]'. 

A significant point which emerges from the above is that we cannot appeal unrestrictedly to set-theoretical reasoning when studying the foundational framework of PA.

Reason: The language of PA has no constant that interprets in any model of PA as the set \textit{\textit{N}} of all natural numbers.

Moreover, the preceding sections show that the Induction Axiom Schema of PA does not allow us to bypass this constraint by introducing an `actual' (or `completed') infinity disguised as an arbitrary constant - usually denoted by \(c\) or \(\infty\) - into either the language, or a putative model, of PA.

\subsection{The case against Goodstein's interpretation of his own set-theoretic reasoning}

The significance of the preceding observations is seen in the following analysis of Goodstein's Theorem, which illustrates why we need to differentiate between theorems in PA and theorems in Ordinal Arithmetic.

Now, in a 1944 paper, R. L. Goodstein\footnote{\cite{Gd44}.}, considers, for any given natural number \(m\), the sequence, \(G(m)\), of natural numbers:
\begin{eqnarray}
G(m) & \equiv & \left\{m_{1<2>}, m_{2<3>}, m_{3<4>}, \ldots \right\} \label{eq:1}
\end{eqnarray}
where \(m_{1<2>}\) denotes the unique hereditary representation of the natural number \(m\) in the natural number base 2:
\begin{quote}
e.\ g.\  \(9_{1<2>} \equiv 1.2^{1.2^{1.2^{0}}+1.2^{0}} + 0.2^{1.2^{1.2^{0}}+0.2^{0}}  + 0.2^{1.2^{0}} + 1.2^{0}\)
\end{quote}
and, for \(n > 2\), \(m_{n<n+1>}\) is defined recursively from \(m_{(n-1)<n>}\) as below.

\subsubsection{The recursive definition of Goodstein's sequence}

For \(n > 2\), we express \(m_{(n-1)<n>}\) syntactically by its hereditary representation as:
\begin{eqnarray}
m_{(n-1)<n>} & \equiv & \sum_{i=0}^{l} a_{i}.n^{i_{<n>}} \label{eq:2}
\end{eqnarray}
where:

\begin{tabbing}
(a) \= \(0 \leq a_{i} < n\) over \(0 \leq i \leq l\); \\
(b) \= \(a_{l} \neq 0\); \\
(c)	\= for each \(0 \leq i \leq l\), the exponent \(i\), too, is expressed syntactically \\ \> by its hereditary representation, \(i_{<n>}\), in the base \(n\); so, also, \\ \> are all of its exponents, and, in turn, all of their exponents, etc.
\end{tabbing}

\addvspace{1ex}
Then, if \(a_{0} \neq 0\):
\begin{eqnarray}
m_{n<n+1>} & \equiv & \sum_{i=1}^{l} (a_{i}.(n+1)^{i_{<n+1>}} ) + (a_{0}-1) \label{eq:3}
\end{eqnarray}
Whilst, if \(a_{i} = 0\) for all \(0 \leq i < k\):
\begin{eqnarray}
m_{n<n+1>} & \equiv & \sum_{i=k+1}^{l} (a_{i}.(n+1)^{i_{<n+1>}} ) + c_{k} \label{eq:4}
\end{eqnarray}
where:
\begin{eqnarray*}
c_{k} & = & a_{k}.(n+1)^{k} - 1 \\
 & = & (a_{k} - 1).(n+1)^{k} + \left\{(n+1)^{k} - 1\right\} \\
 & = & (a_{k}-1).(n+1)^{k} + n\left\{(n+1)^{k-1} + (n+1)^{k-2} \ldots +1\right\}
\end{eqnarray*}
and, so, its hereditary representation, in the base \((n+1)\), is of the form:
\[(a_{k}-1).(n+1)^{k_{<n+1>}} + n\left\{(n+1)^{k_{1<n+1>}} + (n+1)^{k_{2<n+1>}} \ldots +1\right\}\]
where \(k > k_{1} > k_{2} > \ldots \geq 1\).

\subsubsection{Some basic properties of Goodstein's sequence}

For \(n > 1\), we consider the difference:

\begin{quote}
 \(d_{(n-1)} = \left\{m_{(n-1)<n>}- m_{n<n+1>}\right\}\)
\end{quote}

\addvspace{1ex}
Now, if \(a_{0} \neq 0\), we have:
\begin{eqnarray}
d_{(n-1)} & = & \sum_{i=0}^{l} (a_{i}.n^{i_{<n>}}) - \sum_{i=1}^{l} (a_{i}.(n+1)^{i_{<n+1>}} ) - (a_{0}-1) \label{eq:5} \end{eqnarray}

Whilst, if \(a_{i}=0\) for all \(0 \leq i < k\), we have:
\begin{eqnarray}
d_{(n-1)} & = & \sum_{i=k}^{l} (a_{i}.n^{i_{<n>}}) - \sum_{i=(k+1)}^{l} (a_{i}.(n+1)^{i_{<n+1>}} ) - \nonumber \\ & & (a_{k} - 1).(n+1)^{k_{<n+1>}} - \nonumber \\ & & n\left\{(n+1)^{k_{1<n+1>}} + (n+1)^{k_{2<n+1>}} \ldots +1 \right\} \label{eq:6}
\end{eqnarray}
Further, if, in equation (\ref{eq:5}), we replace the base \(<n>\) by the base \(<z>\) in the term:
\begin{eqnarray}
\sum_{i=0}^{l} a_{i}.n^{i_{<n>}} \label{eq:7}
\end{eqnarray}
and the base \(<n+1>\) also by the base \(<z>\) in the term:
\begin{eqnarray}
\sum_{i=k+1}^{l} (a_{i}.(n+1)^{i_{<n+1>}} ) + (a_{0}-1) \label{eq:8}
\end{eqnarray}
we have:
\begin{eqnarray}
d^{\prime}_{(n-1)} & = & \sum_{i=0}^{l} (a_{i}.z^{i_{<z>}}) - \sum_{i=1}^{l} (a_{i}.z^{i_{<z>}}) - (a_{0}-1) \nonumber \\ & = & 1 \label{eq:9}
\end{eqnarray}

Whilst if, in equation (\ref{eq:6}), we replace the bases similarly, we have:
\begin{eqnarray}
d^{\prime}_{(n-1)} & = & \sum_{i=k}^{l} (a_{i}.z^{i_{<z>}}) - \sum_{i=(k+1)}^{l} (a_{i}.z^{i_{<z>}}) - \nonumber \\ & & (a_{k}-1).z^{k_{<z>}} - n\left\{z^{k_{1<z>}} + z^{k_{2<z>}} \ldots +1\right\} \nonumber \\ & = & a_{k}.z^{k_{<z>}} - (a_{k} - 1).z^{k_{<z>}}) - n(z^{k_{1<z>}} + z^{k_{2<z>}} \ldots +1) \nonumber \\ & = & z^{k<z>}-n(z^{k_{1<z>}} + z^{k_{2<z>}} \ldots +1) \label{eq:10}
\end{eqnarray}
where \(k > k_{1} > k_{2} > \ldots \geq 1\).

\addvspace{1ex}
We consider, now, the sequence:
\[Z(m) \equiv (m_{1<2 \mid z>}, m_{2<3 \mid z>}, m_{3<4 \mid z>}, \ldots )\]
obtained from Goodstein's sequence by replacing the base \(<n+1>\), in each of the terms, \(m_{n<n+1>}\), by the base \(<z>\) for all \(n \geq 1\).

\addvspace{1ex}
Clearly, if \(z > n\) for all non-zero terms of the Goodstein sequence, then \( d^{\prime}_{(n-1)} > 0\) in either of the cases - equation (\ref{eq:9}) or equation (\ref{eq:10}).

\addvspace{1ex}
The sequence \(Z(m)\) is, then, a descending sequence of natural numbers, and must terminate finitely.

\addvspace{1ex}
Since \(m_{n<(n+1) \mid z>}\) is also, then, always greater than \(m_{n<n+1>}\), Goodstein's sequence too must terminate finitely in this case.

\addvspace{1ex}
Further, since we can always find a bound \(z > n\) for all non-zero terms of the Goodstein sequence if it terminates finitely, the condition that we can always find some bound \(z > n\) for all non-zero terms of any Goodstein sequence is, thus, equivalent to the assumption that any Goodstein sequence terminates finitely.

\subsubsection{Goodstein's set-theoretic argument}

Now Goodstein's argument is that, if we interpret \(z\) as the first limit ordinal, \( \omega \) and---instead of the natural number sequence, \(Z(m)\)---consider the similarly defined ordinal sequence (where \(1_{o}, 2_{o}, \ldots\) denote the finite ordinals):
\[O_{o}(m_{o}) \equiv \left \{m_{1_{o}<2_{o} \mid \omega>}, m_{2_{o}<3_{o} \mid \omega>}, m_{3_{o}<4_{o} \mid \omega>}, \ldots \right \}\]
then, since we can show that the sequence is monotonically decreasing---and if we accept\footnote{\cite{Me64}, p261.} that there are no infinitely descending sequences of ordinals in any model of ZF\footnote{Note that such an assertion can only be made in a model}---every Goodstein sequence defined similarly as above, but over the finite ordinals:
\begin{eqnarray}
G_{o}(m_{o}) \equiv  \left \{m_{1_{o}<2_{o}>}, m_{2_{o}<3_{o}>}, m_{3_{o}<4_{o}>}, \ldots \right \} \label{eq:1}
\end{eqnarray}
necessarily terminates finitely in every model of ZF.

\subsubsection{Goodstein's conclusion}

Goodstein concludes from this that, since the finite ordinals are isomorphic to the natural numbers, we can assume that every Goodstein sequence defined over the natural numbers must also terminate finitely.

Now, this conclusion implicitly \textit{assumes} that any entailment in ZF involving the ordinal inequality \(\omega >_{o} n_{o}\) for all finite ordinals \(n_{o} \geq 1_{o}\), translates validly as an entailment in PA involving the natural number inequality \(z > n\) for all natural numbers \(n \geq 1\) that index a non-zero term of a Goodstein sequence.

In other words, Goodstein's argument implies that we need not bother to establish a proof that some natural number bound, \(z > n\), always exists for all non-zero terms of any Goodstein sequence in order to conclude that the sequence terminates, even though the latter condition---as we have shown---is, both, necessary and sufficient if a Goodstein sequence defined over the natural numbers is to terminate finitely.

\subsubsection{The flaw in Goodstein's argument}

The flaw in this argument is highlighted if we assume that the Goodstein sequence \(G(k)\) does not terminate for some natural number \(k\).

Then, for any given natural number \(u\), we can find a natural number \(z_{u}\) such that, for \(n < u\), the \(n\)'th term \(k_{n<(n+1) \mid z_{u}>}\) of the sequence \((Z_{u})_{k}\) is always greater than the \(n\)'th term \(k_{n<n+1>}\) of the Goodstein sequence \(G(k)\), and, by definition, there is a 1-1 correspondence between the terms of the two sequences.

Further, it follows from the previous section that these first \(u\) terms of the sequence \((Z_{u})_{k}\) are monotonically decreasing.

Now, by definition, the terms of the monotonically decreasing ordinal sequence in Goodstein's argument:
\[O_{o}(k_{o}) \equiv \left \{k_{1_{o}<2_{o} \mid \omega>}, k_{2_{o}<3_{o} \mid \omega>}, k_{3_{o}<4_{o} \mid \omega>}, \ldots \right \}\]
are also in 1-1 correspondence with the terms of \(G(k)\).

However, by Goodstein's argument, the ordinal sequence \(O_{o}(k_{o})\) must terminate finitely in any model of ZF even in this case, although we have shown that the number of monotonically decreasing terms in the natural number sequence \((Z_{u})_{k}\) can be made arbitrarily large by a suitable choice of \(u\)!

Now, since we have shown that ZF and PA can have no common model, it follows that either ZF is inconsistent or, \textit{even if it is consistent and has a model}, we cannot conclude that the natural number Goodstein sequence \(G(k)\) must terminate finitely simply because Goodstein's corresponding ordinal sequence \(O_{o}(k_{o})\) \textit{would} terminate finitely in a putative model of ZF.

\subsubsection{Analysing Goodstein's Theorem}

Formally, Goodstein's argument would be that---accepting there are no infinitely descending sequences of ordinals in any model of ZF---the descending sequence of ordinal numbers:
\[O_{o}(m_{o}) \equiv \left \{m_{1_{o}<2_{o} \mid \omega>}, m_{2_{o}<3_{o} \mid \omega>}, m_{3_{o}<4_{o} \mid \omega>}, \ldots \right \}\]
must terminate finitely in any model of ZF for any given finite ordinal \(m_{o}\).

\addvspace{1ex}
Hence, the following formula---where \(x_{z<z+1_{o}>}\) is the \(z^{th}\) term of the ordinal Goodstein sequence \(G_{o}(x)\)---is provable in ZF:

\[\left [(\forall x)((x \in \omega) \rightarrow (\exists z)( (z \in \omega) \wedge x_{z<z+1_{o}>} = 0))\right ]\]
since the formula would hold in every model of ZF.

\addvspace{1ex}
Goodstein's Theorem is, then, the conclusion that:
\[(\exists z)(n_{z<z+1>} = 0)\]
holds for any given natural number \(n\) in the `standard' interpretation of PA.

\subsubsection{Why Goodstein's argument is intuitionistically objectionable}

If we accept the formal expression of Goodstein's argument outlined above as faithful, then Goodstein's conclusion, first, mistakenly presumes that the formula `[\((\exists x) \ldots \)]'---which is an abbreviation for the formula `[\(\neg(\forall x)\neg \ldots \)]'---can always be interpreted consistently as Aristotle's particularisation `\((\exists x) \ldots \)'---which is an abbreviation for `There is an/some \(x\) such that \ldots'---under a sound interpretation of any first-order theory that `contains' Peano Arithmetic\footnote{In the sense that every recursive relation can be expressed in the theory in G\"{o}del's sense (cf.\ \cite{Go31} p28(2)).}.

Second, Goodstein's conclusion makes the intuitionistically objectionable presumption that ZF is consistent and has a model, which then allows the inference that the ordinal sequence \(O_{o}(m_{o})\) is well-defined and terminates finitely in some model of ZF.

The intuitionistically objectionable element is highlighted by the following argument.

We can define the ordinal \(w_{o}\) such that it is the smallest ordinal larger than \({m_{i_{o}<(i+1)_{o}>}}\) for all \(i_{o} > 1_{o}\) in the Goodstein sequence defined for the finite ordinal \(m_{o}\):
\begin{eqnarray}
G_{o}(m_{o}) \equiv  \left \{m_{1_{o}<2_{o}>}, m_{2_{o}<3_{o}>}, m_{3_{o}<4_{o}>}, \ldots \right \} \label{eq:1}
\end{eqnarray}

Now, if \(G_{o}(m_{o})\) terminates finitely, then \(w_{o}\) is a finite ordinal. If, however, \(G_{o}(m_{o})\) does not terminate finitely, then it follows that the sequence of ordinals:
\[W_{o}(m_{o}) \equiv \left \{m_{1_{o}<2_{o} \mid w_{o}>}, m_{2_{o}<3_{o} \mid w_{o}>}, m_{3_{o}<4_{o} \mid w_{o}>}, \ldots \right \}\]
does not terminate.

Since \(w_{o}\) is then, by definition, the first limit ordinal \(\omega\)--and so \(W_{o}(m_{o})\) is the monotonically decreasing sequence of ordinals \(O_{o}(m_{o})\)---it follows by the constructive argument given above that whether the ordinal sequence \(O_{o}(m_{o})\) terminates finitely or not in a model of ZF is a consequence of whether the sequence \(G_{o}(m_{o})\) terminates finitely or not in the model.

However, it then follows that Goodstein must appeal to a non-constructive---hence intuitionistically objectionable---interpretation in order to conclude that there is a model of ZF with no infinitely descending sequences of ordinals!

\subsection{ZF is inconsistent}

The blurring of the distinction between theorems in Ordinal Arithmetic and PA---as reflected in the accepted interpretation of Goodstein's argument---can be viewed as another instance of Nathanson's `boss' factor; one, moreover, that may have handicapped the various attempts over the last century to establish that ZF is inconsistent\footnote{For instance, the efforts by the ultra-intuitionist Alexander Yessinin-Volpin's group to show that ZF is inconsistent.}.

I now show how the inconsistency of ZF follows straightforwardly from the Provability Theorem for PA.

\subsubsection{G\"{o}del's \(R(x)\) is a Halting-type of relation}

Prima facie, one reason\footnote{Another is considered in Appendix D.} that \(R(x)\) is constructively computable as always true but not Turing-computable as always true may be that, following G\"{o}del, we define \(R(x)\)---more accurately, \(R(x, p)\)\footnote{G\"{o}del refers to the corresponding formula \([R(x, y)]\) of his formal system P only by its G\"{o}del-number \(q\), to the P-formula \([(\forall x)R(x, y)]\) by its G\"{o}del-number \(p\), and to the P-formula \([R(x, p)]\) by its G\"{o}del-number \(r\).}---as instantiationally equivalent to the primitive recursive relation \(Q(x)\)---more accurately, \(Q(x, p)\)\footnote{G\"{o}del defines \(Q(x, y)\)---with reference to his formal system P and his chosen system of G\"{o}del-numbering (\cite{Go31}, pp.9-14)---as the primitive recursive relation which translates as: `\(x\) is not the G\"{o}del-number of a proof-sequence in P whose last formula is obtained from the formula whose G\"{o}del-number is \(y\) by substituting the P-numeral that represents \(y\) for the P-variable whose G\"{o}del-number is \(19\)' \cite{Go31}, p.24(8.1). He implicitly defines \(R(x, y)\) as the corresponding arithmetical relation by appeal to his Theorem VII \cite{Go31}. It is also implicit in G\"{o}del's reasoning that \([R(x, y)]\) is the P-formula that interprets under the standard interpretation of P as the arithmetical relation \(R(x, y)\). Since G\"{o}del's definition of \(R(x, y)\)---by appeal to his Theorem VII---contains existential assertions expressed by means of the usual existential symbol `\(\exists\)', his reasoning---as also that of Rosser \cite{Ro36}.---\textit{implicitly} assumes that the standard interpretation of P is sound, and so every occurence of the P-symbol \([\exists]\) in the P-formula \([R(x, y)]\) interprets as the concept defined by Hilbert's \(\epsilon\)-symbol---and denoted by the shorthand notation `\(\exists\)'---in the interpreted arithmetical relation \(R(x, y)\).}---in a deliberately self-referential way that involves an implicit circularity.

The circularity is implicit in the definition of \(Q(x, y)\)\footnote{cf.\ \cite{Go31} p24(8.1)}, whose domain includes the G\"{o}del-number of the PA-formula \([R(x, y)]\). Hence, the definition of \(Q(x, p)\)\footnote{The primitive recursive relation denoted by \(Q(x, p)\) translates as: `\(x\) is not the G\"{o}del-number of a proof-sequence in P, the last formula of which is \([(\forall x)R(x, y)]\)'.} implicitly references the range of values of \(R(x, p)\).

As a consequence, there may be some natural number \(n\) such that, when computing \(R(n, p)\), the Turing-machine TM\(_{R}\) that computes \(R(x, p)\) does not halt, but goes into a non-terminating (machine-generated, but not program-generated) loop, on the input \(n\), since an instantaneous tape description\footnote{For a formal definition see, e.g., \cite{Me64}, p.230; Turing (\cite{Tu36}, \S 2) refers to it as ``\ldots the \textit{complete configuration} at that stage" of the TM.} repeats itself.

\begin{quote}
\footnotesize{
\textbf{Non-terminating loop}: A non-terminating loop is defined as any repetition of the instantaneous tape description of a Turing-machine during a computation.
}
\end{quote}

\subsubsection{Turing's Halting problem}

Now, in his seminal paper on computable numbers, Turing considers\footnote{\cite{Tu36}, \S 8.} the Halting problem, which is a more general case of the above, and which can be expressed as the query\footnote{cf.\ \cite{Me64}, p256.}:

\begin{quote}
\textbf{Halting problem for TM}: Given a Turing-machine TM, can one effectively decide, given any instantaneous description alpha, whether or not there is a computation of TM beginning with alpha?
\end{quote}

Turing then shows that the Halting problem is unsolvable by a Turing-machine\footnote{Turing concludes that: ``\ldots there can be no machine \(\mathcal{E}\) which, when supplied with the S.D. of an arbitrary machine \(\mathcal{M}\), will determine whether \(\mathcal{M}\) ever prints a given symbol (\(0\) say)".}. Since a function is Turing-computable if, and only if, it is partial recursive\footnote{\cite{Me64}, p233, Corollary 5.13 \& p237, Corollary 5.15.}, it is essentially unverifiable algorithmically whether or not a Turing-machine that computes a random \(n\)-ary, number-theoretic function will halt on every \(n\)-ary sequence of natural numbers as input, and not go into a non-terminating loop on some input.

Now we note that any Turing-machine can be provided with an effective looping oracle--- in the form of an auxiliary infinite tape\footnote{\cite{Rg87}, p130.}---that will effectively recognise a non-terminating looping situation and abort the computation.

\begin{quote}
\footnotesize{
\textbf{Effective Looping oracle}: An auxiliary device that records every instantaneous tape description at the execution of each machine instruction of a Turing-machine, and compares the current instantaneous tape description with the record. If an instantaneous tape description is repeated, the oracle aborts the computation of the Turing-machine on the impending non-terminating loop, and returns a meta-symbol indicating self-termination as an output of the Turing-machine.
}
\end{quote}

\subsubsection{If ZF is consistent, then the Halting problem is effectively solvable by a Turing-machine}

All of the foregoing reasoning leads to a curious consequence:

\begin{theorem}
: If ZF is consistent, then it can be effectively determined whether, or not, given any partial recursive number-theoretic function \(F(x_{1}, ..., x_{n})\), the Turing-machine TM\(_{F}\) will halt or not on any given natural number sequence \(a_{1}, ...,\) \(a_{n}\) as input.
\end{theorem}

\vspace{+1ex}
\noindent\textbf{Proof}: We assume that \(F\) is obtained from the recursive function \(G\) by means of the unrestricted \(\mu\)-operator\footnote{cf.\ \cite{Me64}, p.121(VI).}, so that \(F(x_{1}, ..., x_{n}) = \mu y(G(x_{1}, ..., x_{n}, y) = 0)\). 

Let \([H(x_{1}, ..., x_{n}, y)]\) express \(\neg (G(x_{1}, ..., x_{n}, y) = 0)\) in PA, so that for any given natural number sequence \(a_{1}, ..., a_{n}, b\):

\vspace{+1ex}
If \(\neg (G(a_{1}, ..., a_{n}, b) = 0)\) is true, then PA proves \([H(a_{1}, ..., a_{n}, b)]\)

If \((G(a_{1}, ..., a_{n}, b) = 0)\) is true, then PA proves \([\neg H(a_{1}, ..., a_{n}, b)]\)

\vspace{+1ex}
We then consider the PA-provability---and truth under a sound interpretation of PA, such as \(\mathcal{I}_{PA(\beta rouwer/Turing)}\)--- of the PA-formula \([H(a_{1}, ..., a_{n}, y)]\) for a given sequence of PA-numerals \([a_{1}], ..., [a_{n}]\). Thus:

\begin{quote}
(a) Let Q1 be the meta-assertion that the Turing-machine TM\(_{G}\) computes \(G(a_{1}, ..., a_{n}, y)\) as \(0\) for some \(y=k\). 

(b) Next, let Q2 be the meta-assertion that the Turing-machine TM\(_{H}\) that computes the arithmetical relation \(H(a_{1}, ..., a_{n}, y)\) returns the symbol for self-termination at the occurrence of a non-terminating loop for some value \(y = k'\).

In other words, given any natural number \(k'\), \(H(a_{1}, ..., a_{n}, k')\) is computable since it is instantiationally equivalent to \(\neg (G(a_{1}, ..., a_{n}, k') = 0)\), but the Turing-machine TM\(_{H}\) cannot always (i.e., for any given \(y\)) compute \(H(a_{1}, ..., a_{n}, y)\) since---as in the case of G\"{o}del's arithmetical relation \(R(x)\) above---there is some value \(y=k'\) for which \(H(a_{1}, ..., a_{n}, k')\) references itself.

(c) Finally, let Q3 be the meta-assertion that TM\(_{G}\) always (i.e., for any given \(y\)) computes \(G(a_{1}, ..., a_{n}, y)\) as non-zero, and that the Turing-machine TM\(_{H}\) also computes \(H(a_{1}, ..., a_{n}, y)\) as always true (i.e., true for any given \(y\)).

Now, if we assume the Provability Theorem for PA, then it follows that \([H(a_{1}, ..., a_{n}, y)]\) is PA-provable. Let \(h\) be the G\"{o}del-number of \([H(a_{1}, ..., a_{n}, y)]\). We consider, then, G\"{o}del's primitive recursive number-theoretic relation \(uB_{PA}v\)\footnote{\cite{Go31}, p.22, def.45.}, which holds if, and only if, \(u\) is the G\"{o}del-number of a proof-sequence in PA for the PA-formula whose G\"{o}del-number is \(v\). It follows that there is some finite \(k''\) such that the Turing-machine TM\(_{B}\) will compute \(uB_{PA}h\) as true for \(u=k''\).
\end{quote}

It follows that Q1, Q2, and Q3 are mutually exclusive and exhaustive, and that one of the parallel trio {TM\(_{G}\), TM\(_{H}\), TM\(_{B}\)} of Turing-machines will always halt.

Now:

\begin{quote}
(i) If TM\(_{G}\) halts on \(y=k\), then TM\(_{F}\) computes \(F(a_{1}, ..., a_{n}) = k\), and halts on the input \(a_{1}, ..., a_{n}\);

(ii) If TM\(_{H}\) halts, then two computable number-theoretic relations that are instantiationally equivalent are not identical. Since PA can be interpreted in ZF by relativisation\footnote{\cite{Me64}, pp.192-193.}, all of the foregoing reasoning can be interpreted in ZF. However, relations are sets in ZF, and two relations that are instantiationally equivalent correspond to two sets that have the same members. Moreover---by a basic ZF axiom---two sets are equal if, and only if, they have the same members. Hence, if TM\(_{H}\) halts, then ZF is inconsistent;

(iii) If TM\(_{B}\) halts, then TM\(_{F}\) does not halt on input \(a_{1}, ..., a_{n}\), and \(F(a_{1}, ..., a_{n})\) is undefined.
\end{quote}

It follows that, if ZF is consistent, then TM\(_{H}\) never halts, and the Halting problem is effectively solvable by a Turing-machine, since TM\(_{F}\) can simulate TM\(_{G}\)+TM\(_{B}\). \(\Box\)

Since the Halting problem is not solvable algorithmically, i.e., by a Turing-machine, we conclude that:

\begin{corollary}
: ZF is inconsistent.
\end{corollary}

\subsubsection{The root of the inconsistency in ZF} 

G\"{o}del constructively defined the arithmetical relation \(R(x)\) in terms of a primitive recursive predicate \(Q(x)\) such that, for any given natural number \(n\), \(R(n) \leftrightarrow Q(n)\). By the ZF axiom of Extensionality, two functions are mathematically (set-theoretically) identical if they are instantiationally equivalent. However, whilst \(Q(x)\) is algorithmically decidable, \(R(x)\) is not. Thus, reflecting Skolem's (apparent) paradox, only one of \(Q(x)\) and \(Q(x)\) defines a set!

\begin{quote}
\footnotesize{
\textbf{Skolem's (apparent) paradox}: In a 1922 address delivered in Helsinki before the Fifth Congress of Scandinavian Mathematicians, Skolem improved upon both the argument and statement of L\"{o}wenheim's 1915 theorem\footnote{\cite{Lo15}, p.235, Theorem 2.}---subsequently labelled as the (downwards) L\"{o}wenheim-Skolem Theorem\footnote{\cite{Sk22}, p.293.}. He then drew attention to a \footnote{\cite{Sk22}, p.295.}:

``\ldots peculiar and apparently paradoxical state of affairs. By virtue of the axioms we can prove the existence of higher cardinalities, of higher number classes, and so forth. How can it be, then, that the entire domain \(B\) can already be enumerated by means of the finite positive integers? The explanation is not difficult to find. In the axiomatization, ``set" does not mean an arbitrarily defined collection; the sets are nothing but objects that are connected with one another through certain relations expressed by the axioms. Hence there is no contradiction at all if a set \(M\) of the domain \(B\) is non-denumerable in the sense of the axiomatization; for this means merely that \textit{within} \(B\) there occurs no one-to-one mapping \(\Phi\) of \(M\) onto \(Z_{o}\) (Zermelo's number sequence). Nevertheless there exists the possibility of numbering all objects in \(B\) , and therefore also the elements of \(M\), by means of the positive integers; of course such an enumeration too is a collection of certain pairs, but this collection is not a ``set" (that is, it does not occur in the domain \(B\))."
}
\end{quote}

\section{Appendix A: Hilbert's interpretation of quantification}

Hilbert interpreted quantification in terms of his \(\varepsilon\)-function as follows\footnote{\cite{Hi27}, p.466.}:

\begin{quote}
IV. The logical \(\varepsilon\)-axiom

13. \(A(a) \rightarrow A(\varepsilon (A))\)

Here \(\varepsilon (A)\) stands for an object of which the proposition \(A(a)\) certainly holds if it holds of any object at all; let us call \(\varepsilon\) the logical \(\varepsilon\)-function. \ldots

1. By means of \(\varepsilon\), ``all" and ``there exists" can be defined, namely, as follows:

    (i) \indent \( (\forall a) A(a) \leftrightarrow A(\varepsilon(\neg A)) \)

    (ii) \( (\exists a) A(a) \leftrightarrow A(\varepsilon(A)) \ldots\)

On the basis of this definition the \(\varepsilon\)-axiom IV(13) yields the logical relations that hold for the universal and the existential quantifier, such as:

   \( (\forall a) A(a) \rightarrow A(b) \) \ldots (Aristotle's dictum),

and:

     \( \neg((\forall a) A(a)) \rightarrow (\exists a)(\neg A(a)) \) \ldots (principle of excluded middle).
\end{quote}

Thus, Hilbert's interpretation of universal quantification --- defined in (i) --- is that the sentence \((\forall x)F(x)\) holds (\textit{under a consistent interpretation \(\mathcal{I}\)}) if, and only if, \(F(a)\) holds whenever \(\neg F(a)\) holds for any given \(a\) (\textit{in \(\mathcal{I}\)}); hence \(\neg F(a)\) does not hold for any \(a\) (\textit{since \(\mathcal{I}\) is consistent}), and so \(F(a)\) holds for any given \(a\) (\textit{in \(\mathcal{I}\)}).

Further, Hilbert's interpretation of existential quantification --- defined in (ii) --- is that \((\exists x)F(x)\) holds (\textit{in \(\mathcal{I}\)}) if, and only if, \(F(a)\) holds for some \(a\) (\textit{in \(\mathcal{I}\)}).

Brouwer's objection to such an unqualified interpretation of the existential quantifier was that, for the interpretation to be considered sound when the domain of the quantifiers under an interpretation is infinite, the decidability of the quantification under the interpretation must be constructively verifiable in some intuitively and mathematically acceptable sense of the term ``constructive"\footnote{\cite{Br08}.}.

Two questions arise:

\begin{tabbing}
\indent (a)	\= Is Brouwer's objection relevant today? \\

\indent (b)	\> If so, can we interpret quantification `constructively'?
\end{tabbing}

\subsection{The standard interpretation \(\mathcal{I}_{PA(Standard/Tarski)}\) of PA}

We consider the standard interpretation, say \(\mathcal{I}_{PA(Standard/Tarski)}\), of first-order Peano Arithmetic (PA).

Now, if [\((\forall x)F(x)\)] and [\((\exists x)F(x)\)] are PA-formulas, and the relation \(F(x)\) is the interpretation under \(\mathcal{I}_{PA(Standard/Tarski)}\) of the PA-formula [\(F(x)\)], then, in current literature:

\begin{tabbing}
\indent (1a) \= [\((\forall x)F(x)\)] is defined as true under \(\mathcal{I}_{PA(Standard/Tarski)}\) if, and only if, \\ \> for any given natural number \(n\), \(F(n)\) holds under \(\mathcal{I}_{PA(Standard/Tarski)}\); \\

\indent (1b) \> [\((\exists x)F(x)\)] is an abbreviation of [\(\neg (\forall x)\neg F(x)\)], and is defined as true \\ \> under \(\mathcal{I}_{PA(Standard/Tarski)}\) if, and only if, it is not the case that, for any \\ \> given natural number \(n\), \(\neg F(n)\) holds under \(\mathcal{I}_{PA(Standard/Tarski)}\); \\

\indent (1c) \> \(F(n)\) holds under \(\mathcal{I}_{PA(Standard/Tarski)}\) for some natural number \(n\) if, and \\ \> only if, it is not the case that, for any given natural number \(n\), \(\neg F(n)\) \\ \> holds under \(\mathcal{I}_{PA(Standard/Tarski)}\).
\end{tabbing}

Since (1a), (1b) and (1c) together interpret [\((\forall x)F(x)\)] and [\((\exists x)F(x)\)] under \(\mathcal{I}_{PA(Standard/Tarski)}\) as intended by Hilbert's \(\varepsilon\)-function, they attract Brouwer's objection.

This answers question (a).

\subsection{A finitary model \(\mathcal{I}_{PA(\beta rouwer/Turing)}\) of PA}

Clearly, the specific target of Brouwer's objection is (1c), which appeals to Platonically non-constructive, rather than intuitively constructive, plausibility.

We can thus re-phrase question (b) more specifically: Can we define an interpretation of PA over [\(\mathcal{N}\)] that does not appeal to (1c)?    

Now, it follows from Turing's seminal 1936 paper on computable numbers that every quantifier-free arithmetical function (\textit{or relation, when interpreted as a Boolean function}) \(F\) defines a Turing-machine \(\textbf{\textit{TM}}_{F}\)\footnote{In the general case, \(\textbf{\textit{TM}}_{F}\) is defined by the quantifier-free expression in the prenex normal form of \(F\).}\footnote{\cite{Tu36}, pp.\ 138-139.} 

Consequently, as shown in a previous section, we can define another interpretation \(\mathcal{I}_{PA(\beta rouwer/Turing)}\) over the structure [\(\mathcal{N}\)], where:

\begin{tabbing}
\indent (2a) \= [\((\forall x)F(x)\)] is defined as true under \(\mathcal{I}_{PA(\beta rouwer/Turing)}\) if, and only if, \\ \> the Turing-machine \(\textbf{\textit{TM}}_{F}\) computes \(F(n)\) as always true (\textit{i.e., as true} \\ \> \textit{for any given natural number} \(n\)) under \(\mathcal{I}_{PA(\beta rouwer/Turing)}\); \\

\indent (2b) \> [\((\exists x)F(x)\)] is an abbreviation of [\(\neg (\forall x)\neg F(x)\)], and is defined as true \\ \> under \(\mathcal{I}_{PA(\beta rouwer/Turing)}\) if, and only if, it is not the case that the \\ \> Turing-machine \(\textbf{\textit{TM}}_{F}\) computes \(F(n)\) as always false (\textit{i.e., as false for} \\ \> \textit{ any given natural number} \(n\)) under \(\mathcal{I}_{PA(\beta rouwer/Turing)}\).
\end{tabbing}

\(\mathcal{I}_{PA(\beta rouwer/Turing)}\) is a finitary model of PA since we have shown that - when interpreted suitably - all theorems of first-order PA are \textit{constructively} true under \(\mathcal{I}_{PA(\beta rouwer/Turing)}\).

This answers question (b).

\subsection{Are \textit{both} interpretations of PA over the structure \(\mathcal{N}\) sound?}

The structure [\(\mathcal{N}\)] can thus be used to define both the standard interpretation \(\mathcal{I}_{PA(Standard/Tarski)}\) and a finitary interpretation \(\mathcal{I}_{PA(\beta rouwer/Turing)}\) for PA.

However, in the finitary model, from the PA-provability of [\(\neg (\forall x)F(x)\)], we may only conclude that \(\textbf{\textit{TM}}_{F}\) does not compute \(F(n)\) as always true under \(\mathcal{I}_{PA(\beta rouwer/Turing)}\).

We may \textit{not} conclude further that \(\textbf{\textit{TM}}_{F}\) must compute \(F(n)\) as false for some natural number \(n\), since \(F(x)\) may be a Halting-type of function that is not Turing-computable\footnote{eg. \cite{Tu36}, pp.\ 132.}. 

In other words, we may not conclude from the PA-provability of [\(\neg (\forall x)F(x)\)] that \(F(n)\) does not hold under \(\mathcal{I}_{PA(\beta rouwer/Turing)}\) for some natural number \(n\).

The question arises: Are both the interpretations \(\mathcal{I}_{PA(Standard/Tarski)}\) and \(\mathcal{I}_{PA(\beta rouwer/Turing)}\) of PA over the structure [\(\mathcal{N}\)] sound?

\subsubsection{PA is \(\omega \)-\textit{in}consistent}

Now, G\"{o}del has constructed\footnote{\cite{Go31}, pp.\ 25(1)} an arithmetical formula, [\(R(x)\)], such that, if PA is assumed simply consistent, then [\(R(n)\)] is PA-provable for any given numeral [\(n\)], but [\((\forall x)R(x)\)] is not PA-provable.

Further, he showed that\footnote{\cite{Go31}, pp.\ 26(2).}, if PA is additionally assumed \(\omega \)-consistent, then [\(\neg (\forall x)R(x)\)] too is not PA-provable.

G\"{o}del defined\footnote{\cite{Go31}, pp.\ 23.} PA as \(\omega\)-consistent if, and only if, there is no PA-formula [\(F(x)\)] for which:

\vspace{+1ex}
(i)	[\(\neg (\forall x)F(x)\)] is PA-provable,

and:

(ii)	[\(F(n)\)] is PA-provable for any given numeral [\(n\)] of PA.

\vspace{+1ex}
However, we have shown that G\"{o}del's reasoning also implies that [\(\neg (\forall x)R(x)\)] is PA-provable, and so PA is \(\omega \)-\textit{in}consistent!

\subsubsection{The interpretation \(\mathcal{I}_{PA(Standard/Tarski)}\) of PA over the structure \(\mathcal{N}\) is \textit{not} sound}

Now, \(R(n)\) holds for any given natural number \(n\), since G\"{o}del has defined \(R(x)\)\footnote{\cite{Go31}, pp.\ 24.} such that \(R(n)\) is instantiationally equivalent to a primitive recursive relation \(Q(n)\) which is computable as true under \(\mathcal{I}_{PA(\beta rouwer/Turing)}\) for any given natural number \(n\) by the Turing-machine \(\textbf{\textit{TM}}_{Q}\).

It follows that we \textit{cannot} admit the standard (\textit{Hilbertian}) interpretion of [\(\neg (\forall x)R(x)\)] under \(\mathcal{I}_{PA(Standard/Tarski)}\) as:

\vspace{+1ex}
\(R(n)\) is false for some natural number \(n\).

\vspace{+1ex}
In other words, the interpretation \(\mathcal{I}_{PA(Standard/Tarski)}\) of PA over the structure [\(\mathcal{N}\)] is \textit{not} sound.

\vspace{+1ex}
However, we \textit{can} interpret [\(\neg (\forall x)R(x)\)] under \(\mathcal{I}_{PA(\beta rouwer/Turing)}\) as:

\vspace{+1ex}
It is not the case that the Turing-machine \(\textbf{\textit{TM}}_{R}\) computes \(R(n)\) as true under \(\mathcal{I}_{PA(\beta rouwer/Turing)}\) for any given natural number \(n\).

\vspace{+1ex}
Moreover, the \(\omega\)-inconsistent PA is consistent with the finitary interpretation of quantification, as in (2a) and (2b) since we have shown above that the interpretation \(\mathcal{I}_{PA(\beta rouwer/Turing)}\) of PA over the structure [\(\mathcal{N}\)] is sound.

\subsubsection{Why the interpretation \(\mathcal{I}_{PA(Standard/Tarski)}\) of PA over \(\mathcal{N}\) is \textit{not} sound}

The reason why the interpretation \(\mathcal{I}_{PA(Standard/Tarski)}\) of PA over the structure [\(\mathcal{N}\)] is \textit{not} sound lies in the fact that, whereas (1b) and (2b) preserve the logical properties of formal PA-negation under interpretation under \(\mathcal{I}_{PA(Standard/Tarski)}\) and \textit{\textbf{B}} respectively, the further \textit{non-constructive} inference in (1c) from (1b) --- to the effect that \(F(n)\) must hold under \(\mathcal{I}_{PA(Standard/Tarski)}\) for some natural number \(n\) --- does not, and is the one objected to by Brouwer\footnote{\cite{Br08}.}.

\section{Appendix B: A Deduction Theorem}

Now, for any first order theory \(\mathcal{K}\), we have the Deduction Theorem:
 
\begin{theorem}
: (Deduction Theorem) If \([A], [B]\) are closed well-formed formulas of \(\mathcal{K}\), and if \(\vdash_{\mathcal{K}}[B]\) when we assume \(\vdash_{\mathcal{K}}[A]\), then \(\vdash_{\mathcal{K}}[A \rightarrow B]\).
\end{theorem}

\noindent \textbf{Proof}: (i) The case \(\vdash_{\mathcal{K}}[B]\) is straightforward, since the \(\mathcal{K}\)-formula [\(A\rightarrow B\)] then interprets as true under any sound interpretation \(\mathcal{I}_{\mathcal{K}}\) of \(\mathcal{K}\). 

\addvspace{+1ex}
(ii) If not \(\vdash_{\mathcal{K}}[B]\) then, if \(\vdash_{\mathcal{K}}[B]\) when we assume \(\vdash_{\mathcal{K}}[A]\), then by definition there is a sequence \([B_{1}], [B_{2}], \ldots, [B_{n}]\), of well-formed \(\mathcal{K}\)-formulas such that \([B_{1}]\) is \([A]\), \([B_{n}]\) is \([B]\) and, for each \(i > 1\), either \([B_{i}]\) is an axiom of \(\mathcal{K}\) or \([B_{i}]\) is a direct consequence by some rules of inference of \(\mathcal{K}\) of the axioms of \(\mathcal{K}\) and some of the preceding well-formed formulas in the sequence. 

\addvspace{+1ex}
(iia) If, now, [\(A\)] is false under a sound interpretation \(\mathcal{I}_{\mathcal{K}}\) of \(\mathcal{K}\), then [\(A\rightarrow B\)] interprets as vacuously true under \(\mathcal{I}_{\mathcal{K}}\). 

\addvspace{+1ex}
(iib) If, however, \([A]\) is true under a sound interpretation \(\mathcal{I}_{\mathcal{K}}\) of \(\mathcal{K}\), then the sequence \([B_{1}], [B_{2}], \ldots, [B_{n}]\) interprets as a deduction squence under \(\mathcal{I}_{\mathcal{K}}\); whence it follows that [\(A\rightarrow B\)] interprets as true under \(\mathcal{I}_{\mathcal{K}}\).

\addvspace{+1ex}
In other words, we cannot have that \([A]\) interprets as true, and \([B]\) as false, under a sound interpretation \(\mathcal{I}_{\mathcal{K}}\) of \(\mathcal{K}\), as this would imply that there is some extension \(\mathcal{K^\prime}\) of \(\mathcal{K}\) in which \(\vdash_{\mathcal{K^\prime}}[A]\), but not \(\vdash_{\mathcal{K^\prime}}[B]\); this would contradict our hypothesis which implies that, in any extension \(\mathcal{K^\prime}\) of \(\mathcal{K}\) in which we have \(\vdash_{\mathcal{K^\prime}}[A]\), the sequence \([B_{1}], [B_{2}], \ldots, [B_{n}]\) yields \(\vdash_{\mathcal{K^\prime}}[B]\).

\addvspace{+1ex}
Hence, [\(A\rightarrow B\)] interprets as true in all models of \(\mathcal{K}\). By a consequence of G\"{o}del's Completeness Theorem for an arbitrary first order theory, it follows that \(\vdash_{\mathcal{K}}[A \rightarrow B]\). \(\Box\)

\subsection{PA is \(\omega\)-\textit{in}consistent}

Now, G\"{o}del has shown\footnote{\cite{Go31}, Theorem VII, pp.29-31.} that we can construct a PA-proposition, \([(\forall x)R(x)]\), such that if the G\"{o}del-number of \([(\forall x)R(x)]\) corresponds to G\"{o}del's \(17Gen \hspace{+.5ex} r\)\footnote{\cite{Go31}, p.25 eqn.(13).}, and if \([(\forall x)R(x)]\) is PA-provable, then there is a PA-numeral \([n]\) such that the PA-formula whose G\"{o}del-number corresponds to G\"{o}del's \(Neg(Sb(r \begin{array}{c} 17 \\ Z(n) \end{array} ))\)---i.e., the PA-formula \([\neg R(n)]\)---is PA-provable if PA is assumed simply consistent\footnote{\cite{Go31}, p.25(1).}.

\vspace{+1ex}
Hence the PA-formula whose G\"{o}del-number corresponds to G\"{o}del's \(Neg(17\) \(Gen \hspace{+.5ex} r)\)---i.e., the PA-formula \([\neg(\forall x)R(x)]\)---is also PA-provable\footnote{Since \(\vdash_{PA}[\neg R(n) \rightarrow (\exists x)\neg R(x)]\)\footnote{\cite{Me64}, p.55 Ex.2(a).}, and \([(\exists x)\neg R(x)]\) is, by definition, \([\neg(\forall x)R(x)]\).} if PA is assumed simply consistent, i.e.:

\addvspace{+1ex}
\(\vdash_{PA}[(\forall x)R(x)] \rightarrow \hspace{+.5ex} \vdash_{PA}[\neg(\forall x)R(x)]\).

\addvspace{+1ex}
By applying the above Deduction Theorem, it follows that:

\addvspace{+1ex}
\(\vdash_{PA}[(\forall x)R(x) \rightarrow \neg(\forall x)R(x)]\)

\addvspace{+1ex}
\(\vdash_{PA}[\neg(\forall x)R(x)]\)

\addvspace{+1ex}
Moreover, we can also prove\footnote{\cite{Go31}, p.26(2).} that if PA is assumed simply consistent, then \(\vdash_{PA}[R(n)]\) for any given natural number \(n\).

Ergo, PA is \(\omega\)-\textit{in}consistent.

\section{Appendix C: Rosser's proposition}

Now, G\"{o}del\footnote{\cite{Go31}, p.24, 8.1.} defines a primitive recursive relation, \( q(x, y) \), that holds if, and only if, \(x\) is the G\"{o}del-number of a well-formed P-formula\footnote{Of his formally defined Peano Arithmetic, P.}, say \( [H(w)] \)---which has a single free variable, \( [w] \)---and \(y\) is the G\"{o}del-number of a P-proof of \( [H(x)] \). So, for any natural numbers \(h, j\):

\addvspace{+1ex}
\noindent (a)	\(q(h, j)\) holds if, and only if, \(j\) is the G\"{o}del-number of a P-proof of \([H(h)]\).

\addvspace{+1ex}
Rosser's argument defines an additional primitive recursive relation, \(s(x, y)\), which holds if, and only if, \(x\) is the G\"{o}del-number of \([H(w)]\), and \(y\) is the G\"{o}del-number of a P-proof of \([\neg H(x)]\). Hence, for any natural numbers \(h, j\):

\addvspace{+1ex}
\noindent (b)	\(s(h, j)\) holds if, and only if, \(j\) is the G\"{o}del-number of a P-proof of \([\neg H(h)]\).

\addvspace{+1ex}
Further, it follows from G\"{o}del's Theorems V\footnote{\cite{Go31}, p.22.} and VII\footnote{\cite{Go31}, p.29.} that the primitive recursive relations \(q(x, y)\) and \(s(x, y)\) are instantiationally equivalent to some arithmetical relations, \(Q(x, y)\) and \(S(x, y)\), such that, for any natural numbers \(h, j\):

\addvspace{+1ex}
\noindent (c)	If \(q(h, j)\) holds, then \([Q(h, j)]\) is P-provable;

\noindent (d)	If \( \neg q(h, j)\) holds, then \([ \neg Q(h, j)]\) is P-provable;

\noindent (e)	If \(s(h, j)\) holds, then \([S(h, j)]\) is P-provable;

\noindent (f)	If \( \neg s(h, j)\) holds, then \([ \neg S(h, j)]\) is P-provable;

\addvspace{+1ex}
Now, whilst G\"{o}del defines \([H(w)]\) as \([( \forall y) \neg Q(w, y)] \), Rosser's argument defines \([H(w)]\) as \([ (\forall y)(Q(w, y) \rightarrow (\exists z)(z \leq y \wedge S(w, z)))] \),

\addvspace{+1ex}
Further, whereas G\"{o}del considers the P-provability of the G\"{o}delian proposition, \( [(\forall y)\neg Q(h, y)]\), Rosser's argument considers the P-provability of the proposition \([(\forall y)(Q(h, y) \rightarrow (\exists z)(z \leq y \wedge S(h, z)))] \). 

\addvspace{+1ex}
We note that, by definition:

\addvspace{+1ex}
\noindent (i)	\(q(h, j)\) holds if, and only if, \(j\) is the G\"{o}del-number of a P-proof of:

\addvspace{+1ex}
\([(\forall y)(Q(h, y) \rightarrow (\exists z)(z \leq y \wedge S(h, z)))]\);

\addvspace{+1ex}
\noindent (ii)	\(s(h, j)\) holds if, and only if, \(j\) is the G\"{o}del-number of a P-proof of:

\addvspace{+1ex}
\([\neg((\forall y)(Q(h, y) \rightarrow (\exists z)(z \leq y \wedge S(h, z))))]\).

\subsection{The formal expression of Rosser's argument}

\noindent (a) We assume, first, that \(r\) is the G\"{o}del-number of some proof sequence in P for the proposition \( [(\forall y)(Q(h, y) \rightarrow (\exists z)(z \leq y \wedge S(h, z)))]\).

\addvspace{+1ex}
\noindent Hence \(q(h, r)\) is true, and \([Q(h, r)]\) is P-provable.

\addvspace{+1ex}
\noindent However, we then have that \([Q(h, r) \rightarrow (\exists z)(z \leq r \wedge S(h, z))]\) is P-provable.

\addvspace{+1ex}
\noindent Further, by Modus Ponens, we have that \([(\exists z)(z \leq r \wedge S(h, z)))]\) is P-provable.

\addvspace{+1ex}
\noindent Now, if P is simply consistent, then \([ \neg ((\forall y)(Q(h, y) \rightarrow (\exists z)(z \leq y \wedge S(h, z))))]\) is not P-provable.

\addvspace{+1ex}
\noindent Hence, \(s(h, n)\) does not hold for any natural number \(n\), and so \( \neg s(h, n)\) holds for every natural number \(n\).

\addvspace{+1ex}
\noindent It follows that \([\neg S(h, n)]\) is P-provable for every P-numeral \([n]\).

\addvspace{+1ex}
\noindent Hence, \([ \neg ((\exists z)(z \leq r \wedge S(h, z)))]\) is also P-provable - a contradiction.

\addvspace{+1ex}
\noindent Hence, \([(\forall y)(Q(h, y) \rightarrow (\exists z)(z \leq y \wedge S(h, z)))]\) is not P-provable if P is simply consistent.

\addvspace{+1ex}
\noindent (b) We assume next that \(r\) is the G\"{o}del-number of some proof-sequence in P for the proposition \([\neg ((\forall y)(Q(h, y) \rightarrow (\exists z)(z \leq y \wedge S(h, z))))]\).

\addvspace{+1ex}
\noindent Hence \(s(h, r)\) holds, and \([S(h, r)]\) is P-provable.

\addvspace{+1ex}
\noindent However, if P is simply consistent, \([(\forall y)(Q(h, y) \rightarrow (\exists z)(z \leq y \wedge S(h, z)))]\) is not P-provable.

\addvspace{+1ex}
\noindent Hence, \(\neg q(h, n)\) holds for every natural number \(n\), and \([\neg Q(h, n)]\) is P-provable for all P-numerals \([n]\). 

\addvspace{+1ex}
\noindent (i)	The foregoing implies \([y \leq r \rightarrow \neg Q(h, y)]\) is P-provable, and we consider \\ \indent the following deduction\footnote{cf.\ \cite{Me64}, p.146.}:

\addvspace{+1ex}
\begin{tabular}{ll}
(1)	\([r \leq k]\)	& \ldots \emph{Hypothesis} \\

(2)	\([S(h, r)]\)	& \ldots \emph{By 3(b)} \\

(3)	\([r \leq k \wedge S(h, r)]\)	& \ldots \emph{From (1), (2)} \\

(4)	\([(\exists z)(z \leq k \wedge S(h, z))]\) 	& \ldots \emph{From (3)}
\end{tabular}

\addvspace{+1ex}
\noindent (ii)	From (1)-(4), by the Deduction Theorem, we have that \([r \leq k \rightarrow (\exists z)(z \leq \\ \indent k \wedge S(h, z))]\) is provable in P for any P-numeral \([k]\);

\addvspace{+1ex}
\noindent (iii)	Now,  \([k \leq r \vee r \leq k]\) is P-provable for any P-numeral \([k]\);

\addvspace{+1ex}
\noindent (iv)	Also, \([(k \leq r \rightarrow \neg Q(h, k)) \wedge (r \leq k \rightarrow (\exists z)(z \leq k \wedge S(h, z)))]\) is P-provable \\ \indent for any P-numeral \([k]\).

\addvspace{+1ex}
\noindent (v)	Hence \([(\neg (k \leq r) \vee \neg Q(h, k)) \wedge (\neg (r \leq k) \vee (\exists z)(z \leq k \wedge S(h, z)))]\) is \\ \indent P-provable for any P-numeral \([k]\).

\addvspace{+1ex}
\noindent (vi)	Hence \([\neg Q(h, k) \vee (\exists z)(z \leq k \wedge S(h, z))]\) is P-provable for any P-numeral \\ \indent \([k]\).

\addvspace{+1ex}
\noindent (vii)	Hence \([(Q(h, k) \rightarrow (\exists z)(z \leq k \wedge S(h, z))]\) is P-provable for any P-numeral \\ \indent \([k]\).

\addvspace{+1ex}
\noindent (viii)	Now, (vii) contradicts our assumption that \([\neg ((\forall y)(Q(h,y) \rightarrow (\exists z)(z \leq \\ \indent y \wedge S(h, z))))]\) is P-provable.

\addvspace{+1ex}
\noindent (ix)	Hence \([\neg ((\forall y)(Q(h, y) \rightarrow (\exists z)(z \leq y \wedge S(h, z))))]\) is not P-provable if \\ \indent P is simply consistent.

\addvspace{+1ex}
However, the claimed contradiction in (viii) \emph{only} follows if we assume that P is \(\omega\)-consistent, and \emph{not} if we assume that P is simply consistent.

\section{Appendix D: A source of true but unprovable arithmetical propositions}

\subsection{G\"{o}del's Theorem V: The Expressibility Theorem}

By G\"{o}del's Theorem V\footnote{\cite{Go31}, p.22.}, every recursive relation \(R(x_{1}, \ldots, x_{n})\) can be expressed in PA by a formula \([F(x_{1}, \ldots, x_{n})]\) such that, for any given \(n\)-tuple of natural numbers \(a_{1}, \ldots, a_{n}\):

\begin{quote}
If \(R(a_{1}, \ldots, a_{n})\) is true, then PA-proves \([F(a_{1}, \ldots, a_{n})]\)

If \(\neg R(a_{1}, \ldots, a_{n})\) is true, then PA-proves \([\neg F(a_{1}, \ldots, a_{n})]\)
\end{quote}

G\"{o}del relies only on the above to conclude---in his Theorem VI\footnote{\cite{Go31}, p.24.}---the existence of an arithmetical proposition that is formally unprovable in a Peano Arithmetic, but true under a sound interpretation of the Arithmetic. 

\subsection{G\"{o}del's Theorem VII: Constructing a true but unprovable arithmetical proposition}

However, we show that it is G\"{o}del's Theorem VII\footnote{\cite{Go31}, p.29.} which, for every recursive relation of the form  \(x_{0}=\phi(x_{1}, \ldots, x_{n})\), provides an actual blueprint for the construction of a PA-formula that is PA-unprovable, but true under any sound interpretation of PA such as \(\mathcal{B}\). 

\subsubsection{Every recursive function is representable in PA}

We note that\footnote{cf.\ \cite{Me64}, pp.131-134.}:

\vspace{+1ex}
\noindent 1. Any recursive function \(f(x_{1}, x_{2})\) can be represented by a PA-formula \([F(x_{1}, x_{2},\) \(x_{3})]\) such that, for any given natural numbers \(k, m, n\), if \(f(k, m) = n\), then:

\begin{quote}
(i) PA proves: \([F(k, m, n)]\)

(ii) PA proves: \([(\exists_{1} x_{3})F(k, m, x_{3})]\)\footnote{The symbol `\([\exists_{1}]\)' denotes uniqueness, in the sense that the PA-formula \([(\exists_{1} x_{3})F(x_{1}, x_{2}, x_{3})]\) is a short-hand notation for the PA-formula \([\neg(\forall x_{3})\neg F(x_{1}, x_{2}, x_{3}) \wedge (\forall y)(\forall z)(F(x_{1}, x_{2}, y) \wedge F(x_{1}, x_{2}, z) \rightarrow y=z)]\).}
\end{quote}

\noindent 2. G\"{o}del's \(\beta\)-function is defined as:

\begin{quote}
\(\beta (x_{1}, x_{2}, x_{3}) = rm(1+(x_{3}+ 1) \star x_{2}, x_{1})\)
\end{quote}

\noindent where \(rm(x_{1}, x_{2})\) denotes the remainder obtained on dividing \(x_{2}\) by \(x_{1}\).

\vspace{+1ex}
\noindent 3. It follows that, for any sequence of values \(f(x_{1}, 0), f(x_{1}, 1), \ldots, f(x_{1}, n)\), we can construct natural numbers \(b, c, i\) such that \(\beta(b, c, i) = f(x_{1}, i)\) for \(0 \leq i \leq n\), where \(c = j\)!, and \(j = max(n, f(x_{1}, 0), f(x_{1}, 1), \ldots, f(x_{1}, n)\).\footnote{cf.\ \cite{Me64}, p.131, Proposition 3.22.}

\vspace{+1ex}
\noindent 4. \(Bt(x_{1}, x_{2}, x_{3}, x_{4})\) is the following representation in PA of \(\beta(x_{1}, x_{2}, x_{3})\)\footnote{cf. \cite{Me64}, p.131.}: 

\begin{quote}
\([(\exists w)(x_{1} = ((1 + (x_{3} + 1)\star x_{2}) \star w + x_{4}) \wedge (x_{4} < 1 + (x_{3} + 1) \star x_{2}))]\).
\end{quote}

\noindent 5. If \(f(x_{1}, x_{2})\) is defined by:

\begin{quote}
(i) \(f(x_{1}, 0) = g(x_{1})\)

(ii) \(f(x_{1}, (x_{2}+1)) = h(x_{1}, x_{2}, f(x_{1}, x_{2}))\)
\end{quote}

\noindent where \(g(x_{1})\) and \(h(x_{1}, x_{2}, x_{3})\) are recursive functions of lower rank\footnote{cf.\ \cite{Me64}, p.132; \cite{Go31}, p.30(2).} that are represented in PA by well-formed formulas \([G(x_{1}, x_{2})]\) and \([H(x_{1}, x_{2}, x_{3}, x_{4})]\), then \(f(x_{1}, x_{2})\) is represented in PA by the following well-formed formula, denoted by \([F(x_{1}, x_{2}, x_{3})]\):

\begin{quote}
\([(\exists u)(\exists v)(((\exists w)(Bt(u, v, 0, w) \wedge G(x_{1}, w))) \wedge Bt(u, v, x_{2}, x_{3}) \wedge (\forall w)(w< x_{2} \rightarrow (\exists y)(\exists z)(Bt(u, v, w, y) \wedge Bt(u, v, (w+1), z) \wedge H(x_{1}, w, y, z)))]\).\footnote{cf.\ \cite{Me64}, p.132.}
\end{quote}

\subsection{Does the well-formed PA-formula \([F(x_{1}, x_{2}, x_{3})]\) represent the recursive function \(f(x_{1}, x_{2})\) strongly in PA?}

The question arises: Is \([(\exists_{1} x_{3})F(x_{1}, x_{2}, x_{3})]\) PA-provable?

\subsubsection{What does ``\((\exists_{1} x_{3})F(k, m, x_{3})\)" assert?}

Now the arithmetical proposition ``\((\exists_{1} x_{3})F(k, m, x_{3})\)" is the assertion that, for any given natural numbers \(k, m\), we can construct natural numbers \(t_{(k, m)}, u_{(k, m)}, v_{(k, m)}\) ---which are functions of \(k, m\)---such that \(\beta(u_{(k, m)}, v_{(k, m)}, 0) = g'(k)\) and, for all \(i<m\), \(\beta(u_{(k, m)}, v_{(k, m)}, i) = h'(k, i, f'(k, i))\), and \(\beta(u_{(k, m)}, v_{(k, m)}, m) = t_{(k, m)}\), where \(f'(x_{1}, x_{2})\), \(g'(x_{1})\) and \(h'(x_{1}, x_{2}, x_{3})\) are any recursive functions that are formally represented by \(F(x_{1}, x_{2},\) \(x_{3}), G(x_{1}, x_{2})\) and \(H(x_{1}, x_{2}, x_{3}, x_{4})\) respectively such that:

\begin{quote}
(i)      \(f'(k, 0) = g'(k)\)

(ii)     \(f'(k, (y+1)) = h'(k, y, f'(k, y))\) for all \(y<m\)

(iii)    \(g'(x_{1})\) and \(h'(x_{1}, x_{2}, x_{3})\) are recursive functions that are assumed to be of lower rank than \(f'(x_{1}, x_{2})\).
\end{quote}

We further note that, for any arbitrarily given sequence of natural numbers \(k_{0}, k_{1}, ... , k_{n}\), we can clearly determine an infinity of values of \(u_{p}, v_{p}, k_{p}\) such that \(\beta(u_{p}, v_{p}, i) = k_{i}\) for all \(0 \leq i \leq n\), and \(\beta(u_{p}, v_{p}, (n+1)) = k_{p}\).

Hence \((\exists_{1} x_{3})F(k, m, x_{3})\) is also the assertion that, for any given natural numbers \(k\) and \(m\), we can always construct some (non-unique) pair of natural numbers \(u_{(k, m)}, v_{(k, m)}\) such that \(\beta(u_{(k, m)}, v_{(k, m)}, i)\) represents the first \(m\) terms, i.e. \(f'(k, 0), f'(k, 1),\) \(\ldots , f'(k, m)\), which are common to every recursive function such as \(f'(x_{1}, x_{2})\) that is formally represented by \([F(x_{1}, x_{2}, x_{3}]\).

We can see that this is constructively provable for any given natural numbers \(k\) and \(m\) since, if \(F(x_{1}, x_{2}, x_{3})\) is a well-defined arithmetical relation\footnote{A critical---and commonly speculated upon---issue that we do not address in this paper is whether the PA-formula \([F(x_{1}, x_{2}, x_{3}]\) can be considered to interpret under a sound interpretation of PA as a well-defined predicate since, for any given \(m\), there are an infinity of functions that are formally represented by \([F(x_{1}, x_{2}, x_{3}]\). This is because the denumerable sequences \(f'(k, 0), f'(k, 1),\) \(\ldots , f'(k, m), m_{p}\)---where \(p>0\) and \(m_{p}\) is not equal to \(m_{q}\) if \(p\) is not equal to \(q\)---are represented by denumerable, distinctly different, functions \(\beta(x_{p_{1}}, x_{p_{2}}, i)\) respectively. There are thus denumerable pairs \(x_{p_{1}}, x_{p_{2}}\) for which \(\beta(x_{p_{1}}, x_{p_{2}}, i)\) represents any given sequence \(f'(k, 0), f'(k, 1),\) \(\ldots , f'(k, m)\).}, it defines the Turing-machine TM\(_{F}\) that can construct the sequence \(f'(k, 0), f'(k, 1), \ldots ,\) \(f'(k, m)\) uniquely and verify the assertion.

\subsubsection{What does ``\((\exists_{1} x_{3})F(x_{1}, x_{2}, x_{3})\)" assert?}

Now, the arithmetical relation ``\((\exists_{1} x_{3})F(x_{1}, x_{2}, x_{3})\)" is the assertion that we can construct natural numbers \(t, u, v\) that are independent of \(x_{1}\) and \(x_{2}\) and such that, for any given \(r, s\), \(\beta(u, v, 0) = g'(r)\) and, for all \(i<s\), \(\beta(u, v, i) = h'(r, i, f'(r, i))\), and \(\beta(u, v, s) = t\).

This, however, is false since---using the above argument---we cannot construct natural numbers \(t, u, v\) that are independent of \(x_{1}\) and \(x_{2}\)\footnote{Since \(u\) is defined as \(j\)!, where \(j = max(n, f(x_{1}, 0), f(x_{1}, 1), \ldots)\), and \(n\) is the number of terms in the sequence \(f(x_{1}, 0), f(x_{1}, 1), \ldots\).}, and such that, for any given \(r, s\), \(\beta(u, v, 0) = g'(r)\) and, for all \(i<s\), \(\beta(u, v, i) = h'(r, i, f'(r, i))\), and \(\beta(u, v, s) = t\).

We conclude that although \((\exists_{1} x_{3})F(k, m, x_{3})\) is true for any given natural numbers \(k, m\), the PA-formula \([(\exists_{1} x_{3})F(k, m, x_{3})]\) is not PA-provable.

We thus have:

\begin{theorem}
Every recursive function is not strongly representable in PA\footnote{cf. \cite{Me64}, p.135, Ex.3.}.
\end{theorem}

\end{document}